\documentclass[final,leqno]{siamltex}
\usepackage{amsmath,amssymb}
\usepackage{graphicx,color}
\usepackage{enumerate}
\usepackage{verbatim}
\usepackage{hyperref}

\newcommand{\om}{\omega}
\newcommand{\Om}{\Omega}
\newcommand{\pa}{\partial}

\newcommand{\bfnu}{\boldsymbol{\nu}}
\newcommand{\bftau}{\boldsymbol{\tau}}

\newcommand{\bfb}{\mathbf{b}}
\newcommand{\bfc}{\mathbf{c}}

\renewcommand{\figurename}{\scriptsize Figure}

\def\be{\begin{equation}}
\def\ee{\end{equation}}
\def\br{\begin{eqnarray}}
\def\er{\end{eqnarray}}

\title{An analysis of sparsity preserving pivot strategies for discontinuous Galerkin methods applied to acoustic scattering}
\author{
Cody Lorton\thanks{Department of Mathematics and Statistics, University of
West Florida, Pensacola, FL 32514, U.S.A.  ({\tt clorton@uwf.edu}).}
\and
Ryan Severance\thanks{Department of Mathematics and Statistics, University of
West Florida, Pensacola, FL 32514, U.S.A.  ({\tt rs104@students.uwf.edu}).}
}

\begin{document}

\maketitle

\begin{abstract}
In this paper we discuss and analyze the sparse structure of matrices associated to the interior penalty discontinuous Galerkin (IP-DG) method applied to the Helmholtz equation.  It is well-known that $LU$-factorization causes fill-in and this paper discusses three pivoting strategies: approximate minimal degree (AMD), nested dissection, and reverse Cuthill-McKee, that can reduce fill-in associated to the $LU$-factorization.  Numerical experiments are included that demonstrate the performance of these pivoting strategies.  These experiments include both uniform and non-uniform mesh structures, the inclusion of a scattering boundary, and both piecewise linear and quadratic solution spaces.
\end{abstract}



\section{Introduction} \label{sec:Intro}

First, we will define the domain used to define the acoustic scattering model 
studied in this paper.  Let $D_1$ and $D_2$ be two open and bounded sets in $
\mathbb{R}^d$ where $d = 2,3$ is the spatial dimension for the problem.  
Further, assume that $D_2 \subset D_1$ such that both $D_1$ and $D_2$ satisfy 
a star-shaped condition involving the same point.  Define the domain $D = D_1/
D_2$ with boundary $\pa D$.  We will decompose the boundary $\pa D$ into $\pa 
D = \Gamma_1 \cup \Gamma_2$ where $\Gamma_1 = \pa D_1$ and $\Gamma_2 = \pa 
D_2$.

Using the domain $D$ defined above we define the following acoustic scattering 
model:
\begin{align}
	- \Delta u - k^2 u &= f, & \mbox{ in } D, \label{eq:PDE1} \\
	\frac{\pa u}{\pa \bfnu} + ik u &= g, & \mbox{ on } \Gamma_1,
	\label{eq:PDE2} \\
	u &= 0, & \mbox{ on } \Gamma_2. \label{eq:PDE3}
\end{align}
\eqref{eq:PDE1} is the acoustic Helmholtz equation and is used to model time-
harmonic acoustic waves in some homogeneous medium.  In this equation 
$u: D \rightarrow \mathbb{C}$ models the displacement of the medium caused by 
the acoustic wave.  The domain $D$ represents the acoustic medium in which the 
wave is traveling. $k$ is the wave number parameter of the wave and is defined 
as $k = \om/c$ where $\om$ is the angular frequency of the wave and $c$ is the 
speed of the wave.  $f,g: D \rightarrow \mathbb{C}$ are given source 
functions.

In the boundary condition \eqref{eq:PDE2} $\bfnu$ refers to the unit outward 
normal vector to the domain $D_1$ and $i = \sqrt{-1}$ is the imaginary unit. 
Typically wave propagation models are posed on unbounded domains with far-
field conditions dictating the behavior of the wave.  These unbounded domains 
are computationally prohibitive.  Thus, the above acoustic scattering model is 
defined on a bounded domain $D$ with boundary $\Gamma_1$ due to the truncation 
of an unbounded domain. To approximate an unbounded domain an absorbing 
boundary condition like \eqref{eq:PDE2} is used \cite{Engquist_Majda_79}.  
\eqref{eq:PDE2} simulates the absorption of the acoustic wave into the 
boundary $\Gamma_1$.

The open and bounded set $D_2$ represents a scattering object in the 
homogeneous media.  Thus, \eqref{eq:PDE3} is a scattering boundary condition.  
In particular, the homogeneous Dirichlet condition \eqref{eq:PDE3} used in 
this paper models a sound-soft scattering object.

This paper is mainly concerned with the case of a large wave number $k$.  It 
is well-known that in this case the Helmholtz differential operator is 
indefinite and this leads to difficulties in analyzing the Helmholtz PDE 
\eqref{eq:PDE1}.  Furthermore, numerical discretization techniques of the 
Helmholtz problem \eqref{eq:PDE1}--\eqref{eq:PDE3} lead to a non-Hermitian 
indefinite system of equations that must be solved to resolve the wave.  This 
leads to subsequent difficulties in the numerical analysis of the 
discretization method as well as difficulties in solving the linear system.

Many discretization methods have been applied to the acoustic Helmholtz 
problem and analyzed.  These include the finite element method (FEM) 
\cite{Aziz_Kellogg_79,Douglas_Santos_Sheen_Bennethum_93,Douglas_Sheen_Santos_94,Melenk_Sauter_11,Wu_13} ,
the plane wave discontinuous Galerkin (PW-DG) method
\cite{Gittelson_Hiptmair_Perugia_09,Hiptmair_Moiola_Perugia_11_2,Moiola_09},
and the ultra weak variational form (UWVF) 
\cite{Buffa_Monk_08,Cessenat_Despres_98,Cessenat_Despres_03,Huttunen_Monk_07,Huttunen_Monk_Collino_Kaipio_04,Luostari_Huttunen_Monk_09}.
It is well-known that one must use a fine enough spatial mesh when defining 
the discretization in order to resolve the wave length $\ell = 2\pi/k$ in each 
coordinate direction.  This leads to a minimum mesh constraint of the form 
$O(kh) = 1$ where $h$ is the mesh size parameter and represents the maximum 
diameter of each element of the partition of the spatial domain. In 
\cite{Ihlenburg_Babuska_95} the authors showed that this minimum constraint is 
necessary for the finite element method (FEM) applied to the 1-dimensional 
acoustic Helmholtz equation.  Furthermore, in \cite{Ihlenburg_Babuska_95} it 
was shown that the $H^1$ error bound for the solution of the FEM applied to 
the 1-d acoustic Helmholtz equation contains a term of the order $k^3h^2$ 
called the pollution term.  This pollution term causes an increase in the 
$H^1$ error as $k$ is increased when $h$ is chosen to satisfy the minimum mesh 
constraint $O(kh) = 1$. The increase in error under the minimum mesh 
constraint is called the pollution effect.   It has been shown that the 
pollution effect is inherent in Helmholtz-type problems and leads to a loss of 
stability in standard discretization techniques 
\cite{Babuska_Sauter_97,Douglas_Santos_Sheen_Bennethum_93}.  Due to this loss 
of stability, a strict mesh constraint of the form $O(k^2h) = 1$ (called the 
asymptotic mesh constraint) 
was required to obtain optimal and quasi-optimal error estimates for the 
acoustic Helmholtz problem 
\cite{Aziz_Kellogg_79,Douglas_Santos_Sheen_Bennethum_93,Douglas_Sheen_Santos_94}.

In \cite{Feng_Wu_09,Feng_Wu_11} an interior penalty discontinuous Galerkin 
(IP-DG) method was developed and analyzed for the acoustic Helmholtz problem 
\eqref{eq:PDE1}--\eqref{eq:PDE3}.  This IP-DG method made use of purely 
imaginary penalty parameters and penalization of the jumps of function values, 
normal derivatives, and tangential derivatives.  In 
\cite{Feng_Wu_09,Feng_Wu_11} this IP-DG method was shown to be unconditionally 
stable. In particular, stability estimates were obtained in both the 
asymptotic and pre-asymptotic mesh regime.  Also, sub-optimal error estimates 
were proven in the pre-asymptotic mesh regime which improve to optimal order 
error estimates in the asymptotic mesh regime.  Numerical experiments in 
\cite{Feng_Wu_09} show that when $k$ is large the IP-DG method outperforms the 
standard FEM in the number of degrees of freedom required to attain a given 
accuracy.  For these reasons, we will focus on this unconditionally stable IP-
DG method in this paper.

To resolve the solution of the acoustic Helmholtz problem 
\eqref{eq:PDE1}--\eqref{eq:PDE3} with $k$ large using standard discretization 
techniques, one must solve a large non-Hermitian, indefinite, and 
ill-conditioned system of linear equations. In \cite{Ernst_Gander_12} it was 
shown that standard iterative solvers applied the the acoustic Helmholtz 
problem do not perform well.  In fact, many do not converge.  Thus, to resolve 
the solution of the acoustic Helmholtz problem a direct linear solver is 
usually employed.  In \cite{Feng_Lin_Lorton_15} an $LU$ direct solver was used 
to obtain an efficient and accurate discretization method for the acoustic 
Helmholtz problem in random media.  

In this paper, we will study the $LU$ decomposition of the system of equations 
obtained from the IP-DG method applied the the acoustic Helmholtz problem 
\eqref{eq:PDE1}--\eqref{eq:PDE3}.  The IP-DG method leads to a linear system 
defined with a sparse global matrix $A$.  It is well-known that the $LU$ 
decomposition leads to ``fill-in" of the sparse system and thus requires more 
CPU memory and leads to a loss of efficiency \cite{George_Ng_85,Gilbert_94}.  
To reduce fill-in, one can use pivoting prior to the $LU$ decomposition.  This 
paper will study three popular sparsity preserving pivoting strategies applied 
to the linear system obtained from the IP-DG method 
\cite{Feng_Wu_09,Feng_Wu_11}: (1) minimum degree pivoting, (2) nested 
dissection, and (3) bandwidth/profile reduction.  We will use numerical 
experiments to compare the fill-in produced by the $LU$ decomposition after 
using these pivoting strategies.

The paper is organized as follows.  Section \ref{sec:IPDG} will detail the IP-
DG method used in this paper.  Section \ref{sec:PivotStrategies} will discuss 
the sparse structure of the global matrix defined by the IP-DG method and 
detail the three sparsity preserving pivoting methods that we will study in 
this paper.  Section \ref{sec:NumericalExp} will be used to discuss multiple 
numerical experiments that were designed to study pivoting to preserve the 
sparsity of the global matrix obtained from the IP-DG method. Section 
\ref{sec:Conclusions} will be used to summarize the results and offer 
conclusions.

\section{Interior penalty discontinuous Galerkin method for acoustic scattering} \label{sec:IPDG}

The goal of this section is to introduce the reader to the interior penalty 
discontinuous Galerkin (IP-DG) method developed and analyzed in 
\cite{Feng_Wu_09,Feng_Wu_11}.  For more in-depth information regarding this 
discretization we encourage the reader to see those papers.  In order to 
define the IP-DG method we will need to introduce standard notation used to 
define DG methods.  For more understanding of discontinuous Galerkin (DG) 
methods we refer the reader to \cite{Riviere_08}.

\subsection{DG Notation} \label{sec:DGNotation}
This subsection will be used to introduce the notation used to define the 
IP-DG method studied in this paper.  Throughout this paper we will use the 
standard complex-valued $L^2$-space norm and inner product notation.  In 
particular, for $\Om \subset D$ and $\Sigma \subset \pa \Om$ we let 
$(\cdot,\cdot)_\Om$ and $\langle \cdot, \cdot \rangle_\Sigma$ denote the 
complex $L^2$-inner product on $\Om$ and $\Sigma$, respectively.  That is,
\begin{align}
	(u,v)_\Om = \int_\Om u \overline{v} dx \qquad  \mbox{and} \qquad 
	\langle u, v\rangle_{\Sigma} = \int_\Sigma u \overline{v} dS. 
	\label{eq:InnerProd}
\end{align}

To discretize the PDE problem \eqref{eq:PDE1}--\eqref{eq:PDE3} using the 
discontinuous Galerkin method, let $\mathcal{T}_h$ be a shape regular 
partition of the domain $D$ parameterized by the size parameter $h$. Typically 
$\mathcal{T}_h$ is a triangulation of a domain in $\mathbb{R}^2$ and $h$ 
specifies the maximum diameter of a triangle $K \in \mathcal{T}_h$.  Example 
triangulations used in later numerical experiments are given throughout 
Section \ref{sec:NumericalExp} (see Figures \ref{fig:Exp1Mesh}, 
\ref{fig:Exp2Mesh},\ref{fig:Exp3Mesh}, and \ref{fig:Exp4Mesh}).

One key characteristic of discontinuous Galerkin methods is the use of an 
energy space $H^1(\mathcal{T}_h) = \prod_{K \in \mathcal{T}_h} H^1(K)$ to 
derive the weak formulation.  This function space includes functions that are 
discontinuous at edge/face boundaries of the partition $\mathcal{T}_h$.  To 
deal with these discontinuities we will need to define special notation that 
is standard in discontinuous Galerkin method.  We begin by defining the 
following sets of edges/faces:
\begin{align*}
	\mathcal{E}^I_h &:= \mbox{ set of all interior edges/faces of } 
	\mathcal{T}_h, \\
	\mathcal{E}^R_h &:= \mbox{ set of all edges/faces of } \mathcal{T}_h 
	\mbox{ on } \Gamma_1, \\
	\mathcal{E}^D_h &:= \mbox{ set of all edges/faces of } \mathcal{T}_h 
	\mbox{ on } \Gamma_2, \\
	\mathcal{E}^{ID}_h &:= \mathcal{E}^I_h \cup \mathcal{E}^D_h.
\end{align*}
For any $e \in \mathcal{E}^I_h$ there exists two cells 
$K,K' \in \mathcal{T}_h$ such that $e = K \cap K'$.  For such an edge/face 
$e = K \cap K'$, let $\bfnu_e$ be the unit normal vector pointing out of $K$ 
if $K$ has larger global label and pointing out of $K'$ if $K'$ has larger 
global label. Also for an edge/face $e = K \cap K'$ we define the jump $[v]$ 
and average $\{ v \}$ operators of a function $v \in H^1(\mathcal{T}_h)$ in 
the following way:
\begin{align*}
	[v] &:= \left\{
	\begin{array}{l} \displaystyle
		v|_K - v|_{K'}, \mbox{ if the global label of } K \mbox{ is larger}, 
		\\
		v|_{K'} - v|_{K}, \mbox{ if the global label of } K' \mbox{ is 
		larger},
	\end{array}
	\right.\\
	\{ v \} & := \frac{1}{2}\left( v|_{K} + v|_{K'} \right).
\end{align*}
For $e \in \mathcal{E}^R_h \cup \mathcal{E}^D_h$ we define $\bfnu_e$ to be the 
unit outward normal vector on $\pa D$ and $[v] = \{v \} := v|_e$.  For 
$e \in \mathcal{E}^I_h \cup \mathcal{E}^R_h \cup \mathcal{E}^D_h$, let 
$\{ \bftau^j_e\}_{j=1}^{d-1}$ be an independent set of unit tangential vectors 
to $e$ and $h_e = \rm{diam}(e)$.

To define the IP-DG method for \eqref{eq:PDE1}--\eqref{eq:PDE3}, we define the 
finite dimensional solution space 
$V^p_h := \prod_{K \in \mathcal{T}_h} P_p(K)$, where $p$ is a positive 
integer.  In other words $V^p_h$ is the set of piecewise polynomials of degree 
$p$ over the partition $\mathcal{T}_h$.  In numerical experiments later in 
this paper we will focus on the cases of $p = 1, 2$, i.e. solutions that are 
piecewise linear or quadratic over the partition $\mathcal{T}_h$.  Over the 
space $V^p_h \times V^p_h$ we define the following sesquilinear form from \
\cite{Feng_Wu_09,Feng_Wu_11}:
\begin{align}
	\label{eq:SesqForm} a_h(u,v) &:= b_h(u,v) - k^2(u,v)_D + ik \langle u, v 
	\rangle_{\Gamma_1} \\
	& \qquad + i\big(J_0(u,v) + J_1(u,v) + L_1(u,v)\big) \qquad \forall
	u,v \in V^p_h, \notag
\end{align}
where
\begin{align*}
	b_h(u,v) &:= \sum_{K \in \mathcal{T}_h} (\nabla u, \nabla v)_K - 
	\sum_{e \in \mathcal{E}^{ID}_h} \left( \left\langle \left\{ \frac{\pa u}
	{\pa \bfnu_e} \right\}, [v] \right\rangle_e + 
	\left\langle [u], \left\{ \frac{\pa v}{\pa \bfnu_e} \right\} 
	\right\rangle_e \right), \\
	J_0(u,v) &:= \sum_{e \in \mathcal{E}^{ID}_h} \frac{\gamma_{0,e}}{h_e} 
	\big\langle [u],[v] \big\rangle, \\
	J_1(u,v) &:= \sum_{e \in \mathcal{E}^{I}_h} \gamma_{1,e} h_e 
	\left\langle \left[ \frac{\pa u}
	{\pa \bfnu_e} \right],\left[ \frac{\pa v}
	{\pa \bfnu_e} \right] \right\rangle, \\
	L_1(u,v) &:= \sum_{e \in \mathcal{E}^{ID}_h}
	\sum_{j=1}^{d-1} \frac{\beta_{1,e}}{h_e}
	\left\langle \left[ \frac{\pa u}
	{\pa \bftau^j_e} \right],\left[ \frac{\pa v}
	{\pa \bftau^j_e} \right] \right\rangle
\end{align*}

The IP-DG method developed and analyzed in \cite{Feng_Wu_09,Feng_Wu_11} and 
used in this paper is defined as: Find $u \in V^p_h$ such that 
\begin{align}
	a_h(u,v) = (f,v)_D + \langle g,v \rangle_{\Gamma_1}, 
	\qquad \forall v \in V^p_h. \label{eq:IPDG}
\end{align}
The terms $J_0,J_1,L_1$ are referred to in the DG literature as penalty terms 
with $\gamma_{0,e},\gamma_{1,e},\beta_{1,e}$ called penalty parameters.  This 
IP-DG method is characterized by penalization of jump discontinuities in the 
function values, normal derivatives, and tangential derivatives at cell 
boundaries in the partition $\mathcal{T}_h$.  The penalty parameters 
$\gamma_{0,e},\gamma_{1,e},\beta_{1,e}$ can be complex values with positive 
real part, but it has been shown that the imaginary part of these parameters 
do not improve performance \cite{Feng_Wu_09}.  Thus, we will assume the 
penalty parameters $\gamma_{0,e},\gamma_{1,e},\beta_{1,e}$ are positive real 
numbers. With this in mind, the penalty terms in \eqref{eq:IPDG} characterize 
purely imaginary penalization which is a unique characteristic of this method.

The key feature of this method, is the unconditional stability of the method.  
In \cite{Feng_Wu_09, Feng_Wu_11}, stability estimates were proven for any 
positive parameters $k,h,\gamma_{0,e},\gamma_{1,e},\beta_{1,e}$.  In contrast,  
for standard discretization methods such as the finite element method, 
stability estimates are only proven in the asymptotic mesh regime 
$O(k^2 h) = 1$.  The unconditional stability of the IP-DG method in 
\eqref{eq:IPDG} leads to unconditional unique solvability of the IP-DG problem 
$\eqref{eq:IPDG}$.  The IP-DG method is also proven to be optimally convergent 
in the asymptotic mesh regime and sub-optimally convergent in the 
preasymptotic mesh regime \cite{Feng_Wu_09,Feng_Wu_11}.

\section{Sparsity preserving pivot strategies} \label{sec:PivotStrategies}
It is well-known that the IP-DG method \eqref{eq:IPDG} is equivalent to 
solving a system of linear equations of the form
\begin{align}
	A \bfc = \bfb. \label{eq:LinSystem}
\end{align}
In this section, we will discuss the sparse properties of this linear system 
as well as sparsity preserving pivoting strategies that can be used to reduce 
fill-in of the $LU$ decomposition.  

First, we note that the IP-DG function space $V^p_h$ is a finite dimensional 
subspace of $H^1(\mathcal{T}_h)$, thus there exists a basis 
$\{ \phi_j \}_{j=1}^N$ of $V^p_h$.  Recall, that $V_h^1$ is the space of 
piecewise linear polynomials over the partition $\mathcal{T}_h$.  A typical 
basis function used in IP-DG methods for $\Om \subset \mathbb{R}^2$ takes the 
form
\begin{align}
	\phi_j = \left\{
		\begin{array}{l r} \displaystyle
			a x + b y + c & \mbox{ in } K^*, \\
			0 & \mbox{ in } K \neq K^*,
		\end{array}
	\right. \label{eq:BasisEx}
\end{align}
where $K^*$ is a single cell in the partition $\mathcal{T}_h$.
In other words, a typical basis function is a polynomial on a single cell that 
is 0 on all other cells of $\mathcal{T}_h$.

Given a basis $\{ \phi_j \}_{j=1}^N$ of $V^p_h$ the linear system in 
\eqref{eq:LinSystem} is defined by 
$A \in \mathbb{C}^{N \times N}$, $\bfc, \bfb \in \mathbb{C}^N$ where
\begin{align*}
	[A]_{i,j} = a_h(\phi_j, \phi_i), \qquad b_i = (f,\phi_i)_D + \langle g,
	\phi_i \rangle_{\Gamma_1},
\end{align*}
and the IP-DG solution $u \in V^p_h$ is defined as
\begin{align*}
	u = \sum_{i = 1}^N c_i \phi_i.
\end{align*}
Thus, for a basis comprised of functions $\phi_i$ with support on a single cell $K \in \mathcal{T}_h$,
\begin{align*}
	[A]_{i,j} = a_h(\phi_j, \phi_j) \neq 0,
\end{align*}
if and only if $\phi_i$ and $\phi_j$ are basis functions associated to the 
same cell or neighboring cells of the partition $\mathcal{T}_h$.  Thus, the 
matrix $A$ is a sparse matrix.  Also, since 
$a_h(\phi_j,\phi_i) = a_h(\phi_i,\phi_j)$ for real-valued basis functions in 
\eqref{eq:SesqForm}, then $A$ is symmetric.  The distribution of nonzero 
element in $A$ depends on the enumeration of the partition $\mathcal{T}_h$ as 
well as the basis functions used.  

As stated earlier in the paper, we will focus on an $LU$ factorization as the 
method for solving the system in \eqref{eq:LinSystem}.  It is well known that 
the $LU$ factorization causes fill-in of sparse systems.  To mitigate this 
fill-in, we will use different permutation matrices $P$ such that the permuted 
matrix $PAP^T$ has an $LU$ factorization with less fill-in. In this paper, 
will test three sparsity preserving pivoting strategies for the system 
\eqref{eq:LinSystem}.  These are minimum degree pivoting (MDP), nested 
dissection (ND), and bandwidth/profile reducing pivoting.  These methods were 
developed for the Cholesky factorization of a symmetric positive definite 
matrix, but can also be applied to the $LU$ factorization of a more general 
square matrix.  \cite{Davis_2006} discusses these pivoting strategies and 
others in depth.

MDP was developed as a sparsity preserving pivoting strategy in calculating a 
right-looking sparse Cholesky factorization.  In particular, MDP is a greedy 
algorithm with the goal of choosing the sparsest pivot row and column in each 
step of the Cholesky factorization.  MDP was first introduced by Markowitz in 
\cite{Markowitz_57}.  In this paper we use the approximate minimal degree 
(AMD) pivoting algorithm developed by Amestoy, Davis, and Duff 
\cite{Amestoy_Davis_Duff_96,Amestoy_Davis_Duff_04}.  For numerical experiments 
discussed in Section \ref{sec:NumericalExp} we made use of Matlab's 
\texttt{amd} function \cite{Matlab_AMD} to obtain the AMD permutation matrix 
$P$.

The nested dissection (ND) permutation was developed by George 
\cite{George_73}.  The ND permutation was designed specifically to reduce 
fill-in due to Cholesky factorization applied to linear systems generated by 
the finite element method.  ND might be a good choice for the matrices 
discussed in this paper due to the similarities between the matrix generated 
by the finite element method and the IP-DG method.  Duff, Erisman, and Reid 
\cite{Duff_Erisman_Reid_76} made an early comparison of MDP and ND 
permutations.  For numerical experiments discussed in Section 
\ref{sec:NumericalExp} we made use of Matlab's \texttt{dissect} function 
\cite{Matlab_ND} to obtain the ND permutation matrix $P$.

The last class of pivoting strategies that we consider is one in which $P$ is 
generated such that the bandwidth/profile of $PAP^T$ is decreased, thus 
reduces the fill-in from the $LU$ decomposition.  One of the first pivoting 
techniques developed for bandwidth/profile reduction was developed by Cuthill 
and McKee in \cite{Cuthill_McKee_69}.  Later Liu and Sherman introduced the 
reverse Cuthill McKee (RCM) method \cite{Liu_Sherman_76}.  RCM reverses the 
Cuthill-McKee ordering which can further reduce the profile of the matrix.  
Chan and George developed an efficient RCM algorithm in \cite{Chan_George_80}.  
For numerical experiments discussed in Section \ref{sec:NumericalExp} we made 
use of Matlab's \texttt{symrcm} function \cite{Matlab_RCM} to obtain the RCM 
permutation matrix $P$.  

\section{Numerical experiments} \label{sec:NumericalExp}
In this section, we present a number of numerical experiments with the intent 
to compare the performance of approximate minimal degree (AMD), nested 
dissection (ND), and reverse Cuthill-McKee (RCM) pivoting strategies.  In 
particular, we will observe the percent of fill-in when an $LU$ factorization 
is used without these pivoting strategies and after these pivoting strategies 
have been used.  The results will be presented using tables that specify the 
total number of nonzero entries compared to the total number of entries, as 
well as a table that gives the percent of fill-in.  In addition to these 
tables, plots will be used to demonstrate the sparsity structure of the matrix 
$A$ before and after pivoting as well as the $LU$ decomposition before and 
after pivoting.  These plots were generated using Matlab's \texttt{spy} 
function.  To save space, the plots of the $LU$ factorization show the 
combined $LU$ factorization defined as $L + U - I$ where $I$ is the identity 
matrix of the same size as $A$ and $A = LU$.

In all experiments in this section we let $D_1$ be the unit square at the 
origin $D_1 = [-0.5,0.5]\times[-0.5,0.5]$ and $k = 5$.  Since our intent is 
only to study the matrix $A$ defined in \eqref{eq:LinSystem} and its $LU$ 
factorization it is unnecessary to specify the functions $f,g$ in 
\eqref{eq:PDE1} and \eqref{eq:PDE2}.   Since the structure of the matrix $A$ 
depends on the structure of the mesh $\mathcal{T}_h$ and the function space 
$V^p_h$ used in the IP-DG method \eqref{eq:IPDG}, our experiments are conducted 
on using different mesh structures and function spaces.  Following 
\cite{Feng_Wu_09}, we chose the following mesh parameters for each experiment
\begin{align*}
	\gamma_{1,e} = 0.1, \qquad \gamma_{0,e} = (k^2h_e)^{2/3} 
	\gamma_{1,e}^{1/3}, \qquad \beta_{1,e} = 1.
\end{align*}

To produce the mesh $\mathcal{T}_h$ and matrix $A$ associated to the IP-DG 
method, FreeFEM++ \cite{Hecht_12} was used.  All matrix analysis was done 
using Matlab's built-in functions.  All experiments were conducted on the same 
Mac computer with a 2 GHz Intel Core i7 processor and 8 GB 1600 MHz DDR3 RAM.

\subsection{Numerical experiment 1} \label{sec:NumericalExp1}
\begin{figure}[!htb]
\includegraphics[width = .5\textwidth]{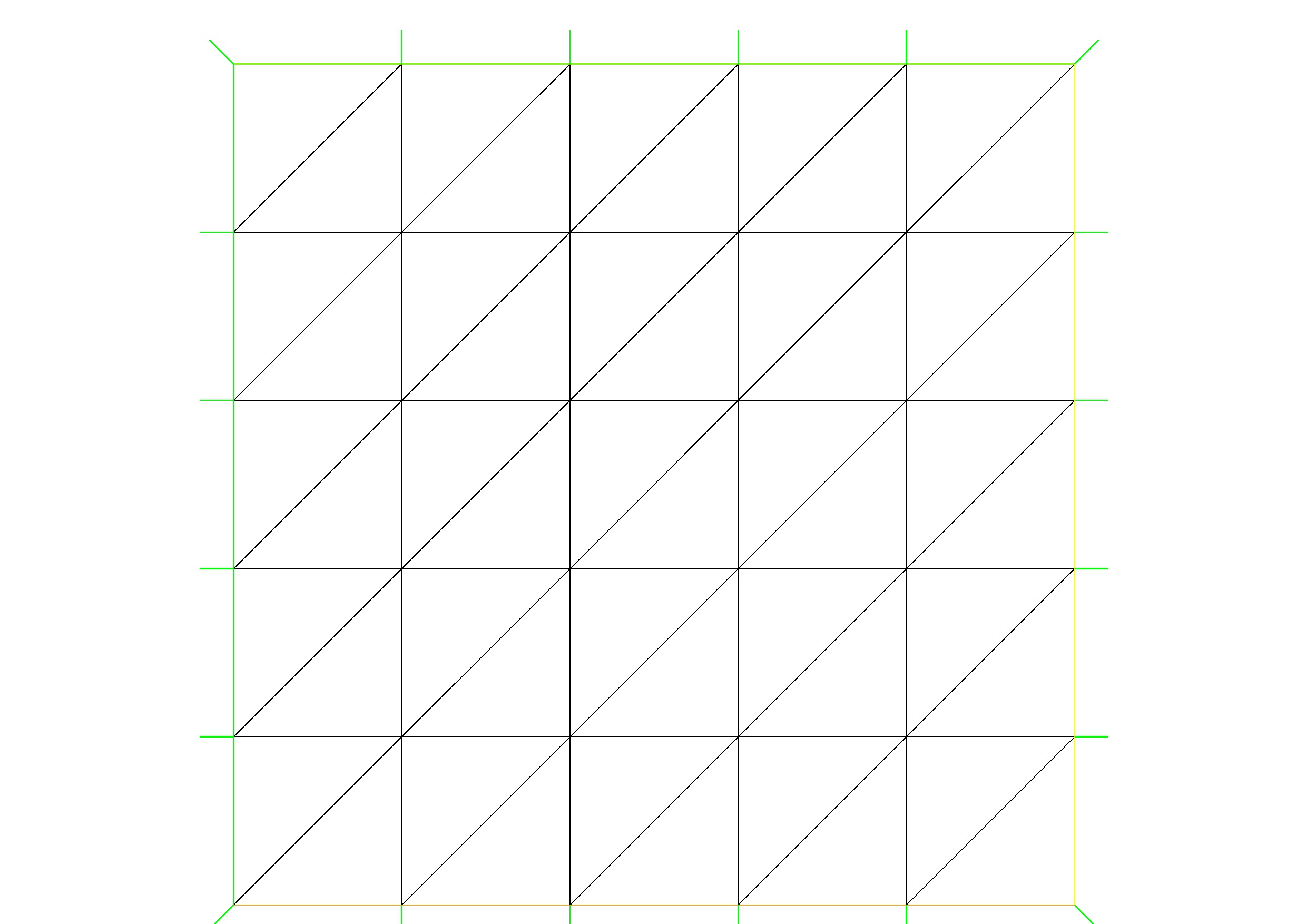}
\includegraphics[width = .5\textwidth]{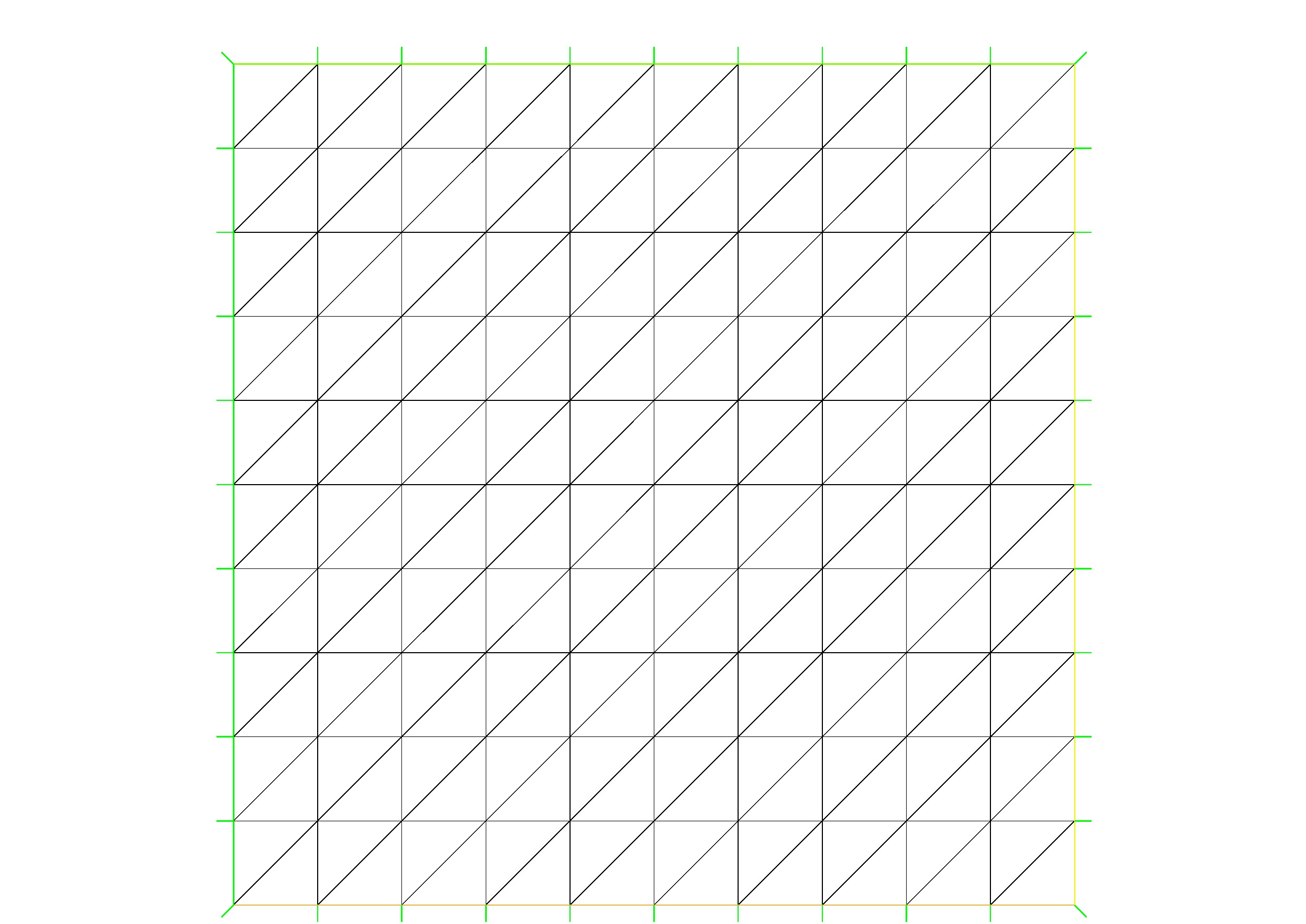}
\caption{\label{fig:Exp1Mesh} Mesh used in experiment 1 with $n = 5$ (left) 
and $n = 10$ (right).}
\end{figure}

In this first set of experiments the scattering portion of the domain was 
chosen to be $D_2 = \emptyset$.  Thus, $D = D_1$ and $\Gamma_2 = \emptyset$.  
To define the mesh $\mathcal{T}_h$ we used a uniform triangulation of $D$ with 
$n = 5, 10, 15, 20$ intervals along each side of $D$.  The mesh size parameter 
is then defined as $h = 1/n$.  Examples of the mesh used in this experiment 
are displayed in Figure \ref{fig:Exp1Mesh}.  In this set of experiments the 
IP-DG function space $V^1_h$ was used, i.e. the set of piecewise linear 
functions across the triangulation $\mathcal{T}_h$.

\begin{figure}[!htb]
\includegraphics[width = .5\textwidth]{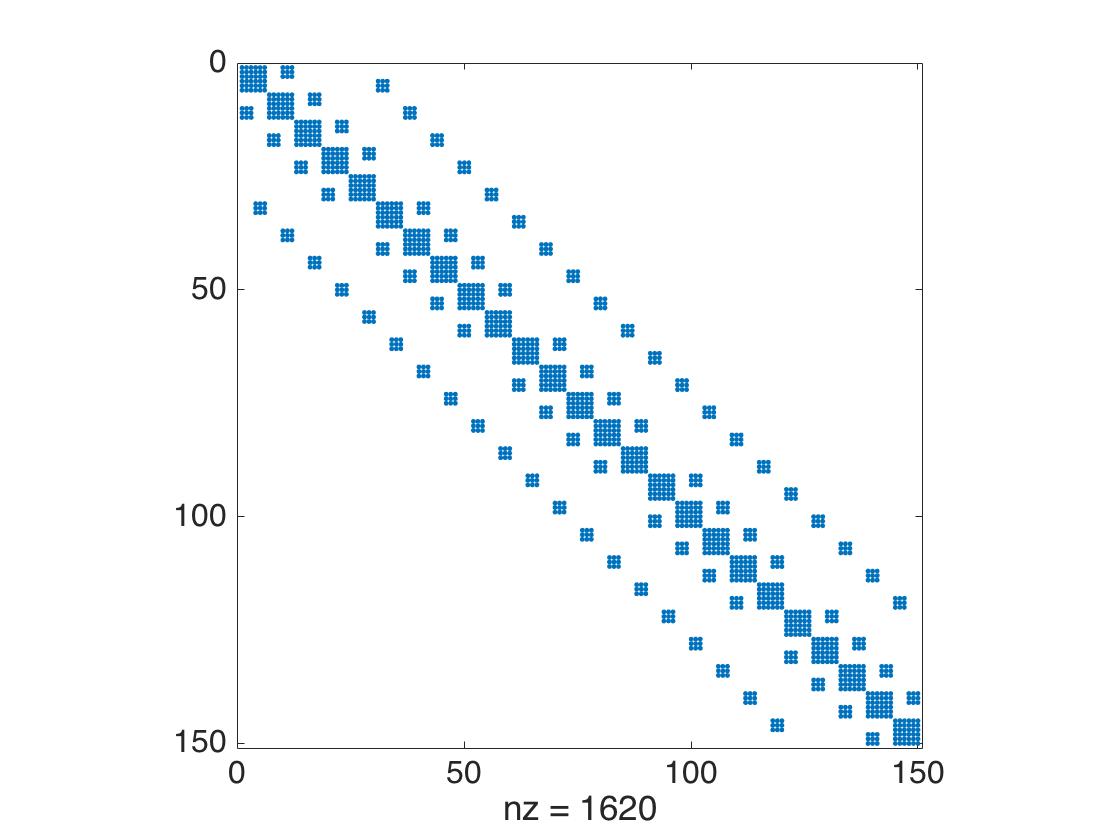}
\includegraphics[width = .5\textwidth]{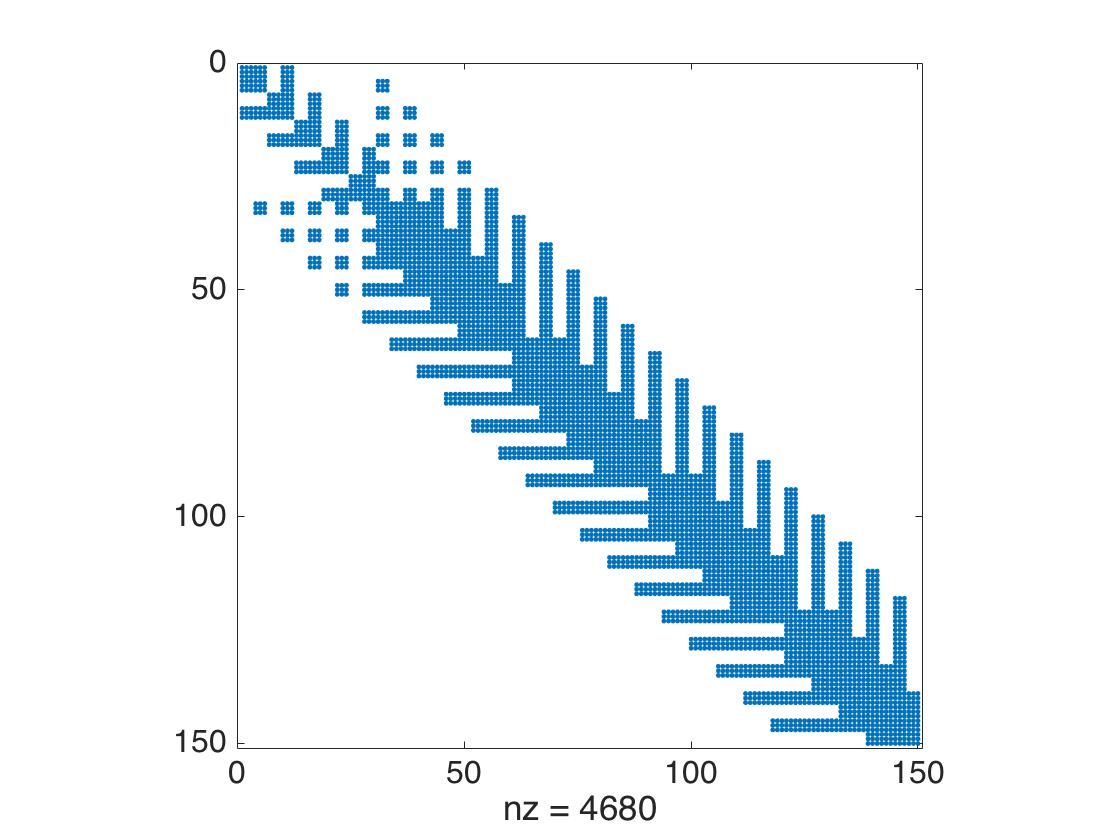}
\caption{\label{fig:Exp1NoPivot5} (Left) sparsity structure of the global 
matrix A produced by the IP-DG method with $n = 5$. (Right) sparsity structure 
of the combined LU decomposition of A.}
\end{figure}

\begin{figure}[!htb]
\includegraphics[width = .5\textwidth]{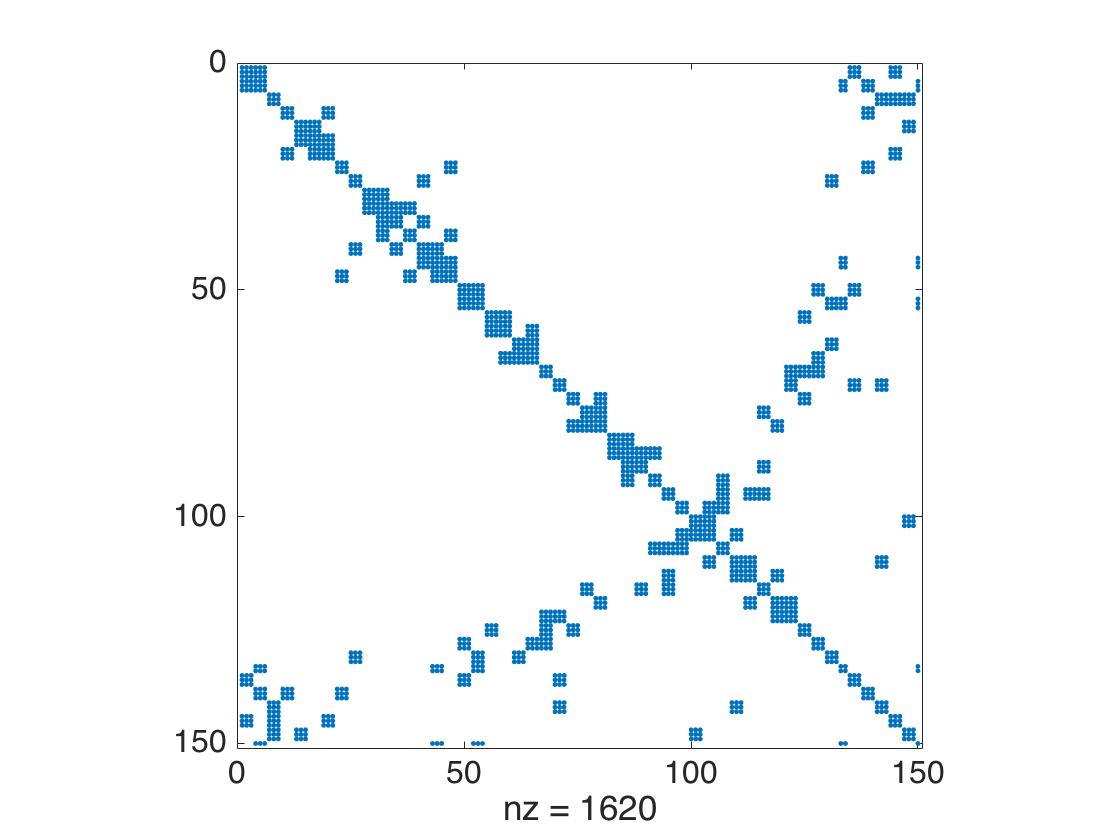}
\includegraphics[width = .5\textwidth]{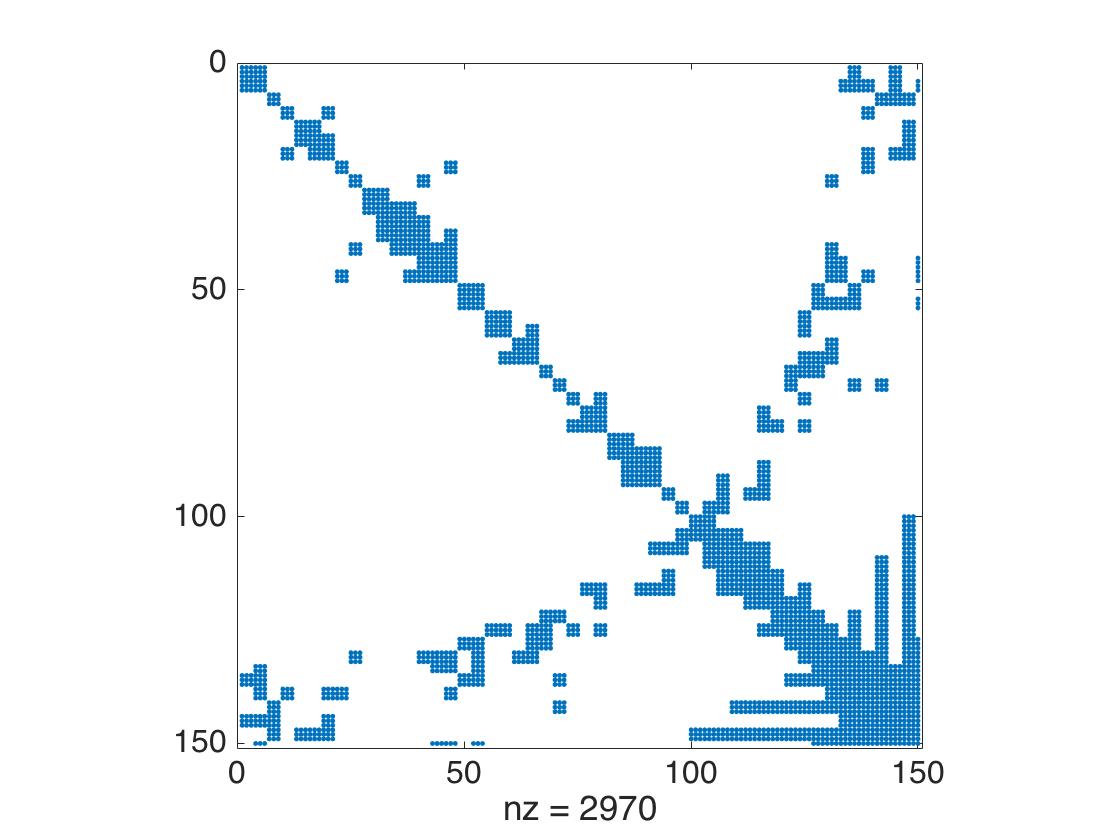}
\caption{\label{fig:Exp1AMD5} (Left) sparsity structure of the global matrix 
after AMD pivoting. (Right) sparsity structure of the combined LU 
decomposition the global matrix after AMD pivoting.}
\end{figure}

\begin{figure}[!htb]
\includegraphics[width = .5\textwidth]{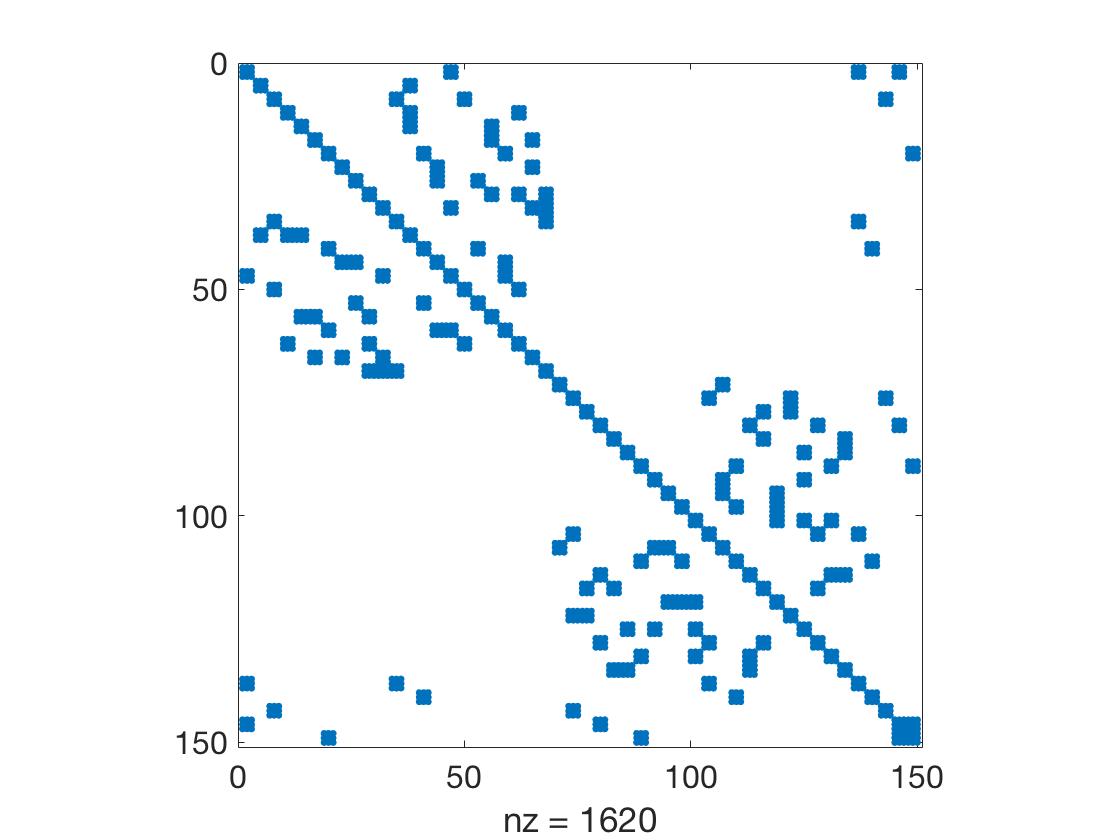}
\includegraphics[width = .5\textwidth]{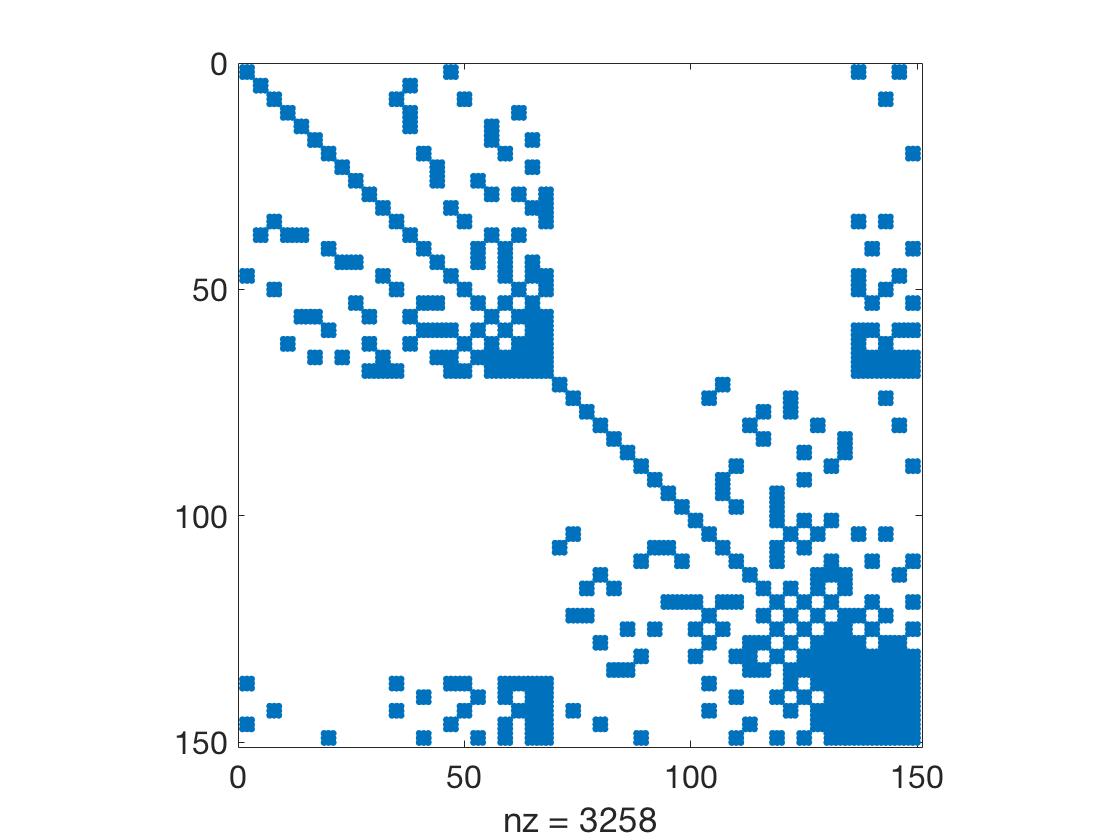}
\caption{\label{fig:Exp1ND5} (Left) sparsity structure of the global matrix 
after ND pivoting. (Right) sparsity structure of the combined LU decomposition 
the global matrix after ND pivoting.}
\end{figure}

\begin{figure}[!htb]
\includegraphics[width = .5\textwidth]{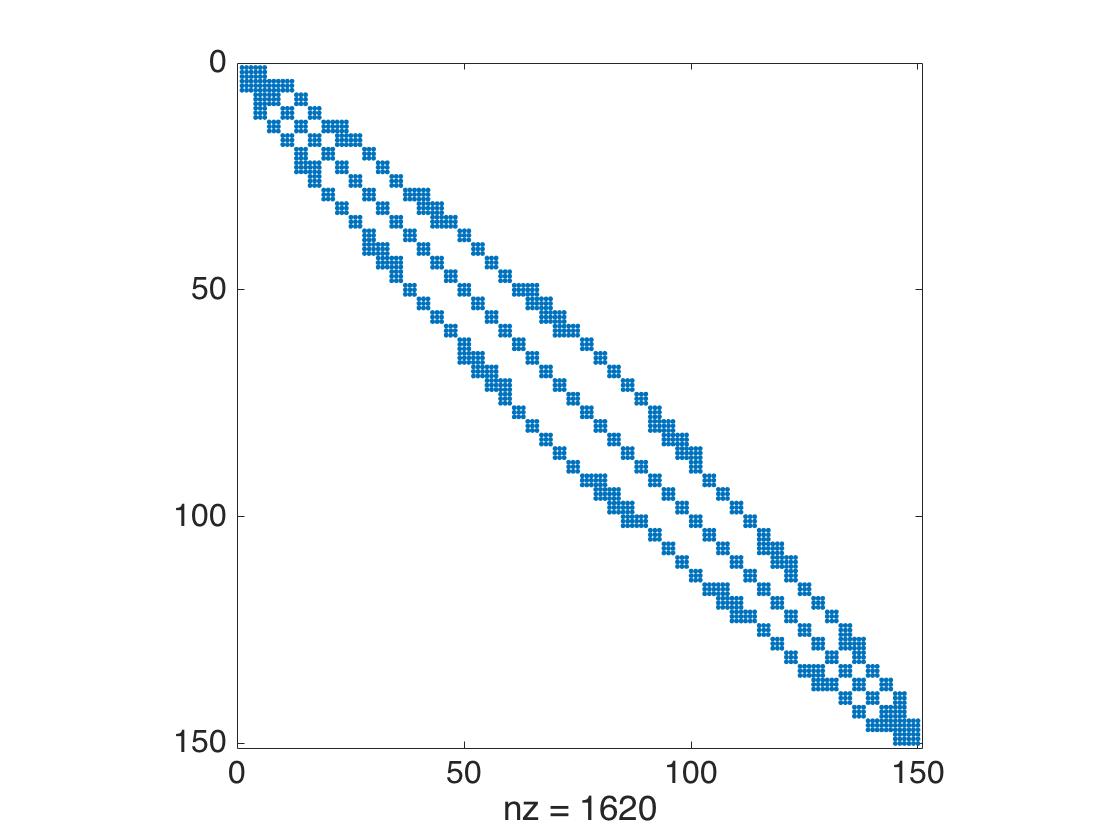}
\includegraphics[width = .5\textwidth]{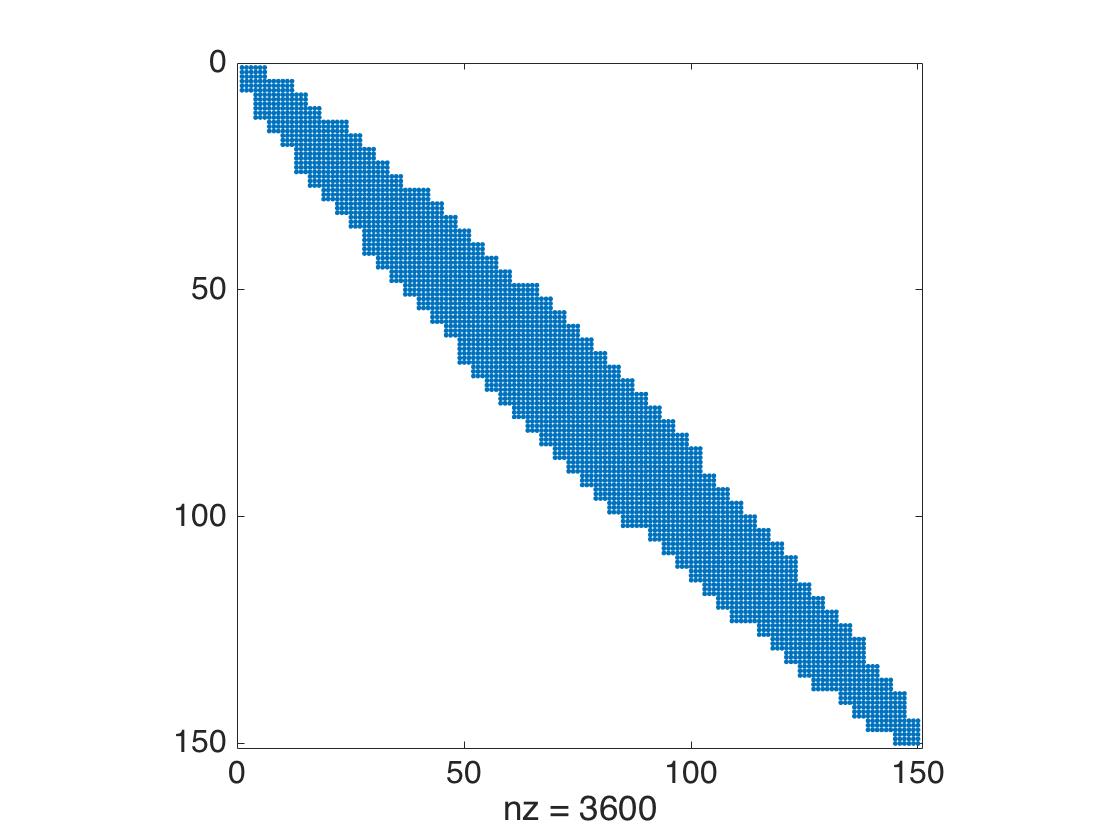}
\caption{\label{fig:Exp1RCM5} (Left) sparsity structure of the global matrix 
after RCM pivoting. (Right) sparsity structure of the combined LU 
decomposition the global matrix after RCM pivoting.}
\end{figure}

Figure \ref{fig:Exp1NoPivot5} demonstrates the sparse structure of the matrix 
$A$ along with its combined $LU$ decomposition for $n = 5$. At the bottom of 
each figure the number of nonzero entries is given as $nz$. This figure 
demonstrates the fill-in associated to the $LU$ factorization.  Figures 
\ref{fig:Exp1AMD5}--\ref{fig:Exp1RCM5} show the sparsity structure of the 
permutation $PAP^T$ produced along with the $LU$ factorization after 
permutation using AMD, ND, and RCM, respectively.  From these images it is 
clear that the pivoting strategies give different sparsity structure and 
varying results with respect to fill-in reduction. Recall that RCM is designed 
to the reduce bandwidth/profile of the matrix $A$.  We can see this bandwidth/
profile reduction clearly in Figure \ref{fig:Exp1RCM5}.

Tables \ref{tab:Exp1_a} and \ref{tab:Exp1_b} summarize the results of 
experiment 1.  In particular, Table \ref{tab:Exp1_a} gives the number of total 
entries of the matrix $A$, the number of nonzero entries of $A$, the number of 
nonzero entries in the $LU$ factorization of $A$, and the number of nonzero 
entries of $LU$ factorization of $PAP^T$ after AMD, ND, and RCM pivoting.  
Table \ref{tab:Exp1_b} presents the same data from Table \ref{tab:Exp1_a} as 
percentage of fill-in.  From Table \ref{tab:Exp1_b} it is clear that $LU$ 
factorization produces a large amount of fill-in.  The percent of fill-in 
decreases as the mesh is refined.  It is also clear that each pivoting 
strategy decreases the percent of fill-in which leads to improved performance.  
In this experiment both AMD and ND pivoting decreases fill-in better than RCM 
pivoting.  For a coarse mesh $n = 5$ AMD performs slightly better than ND, but 
as the mesh is refined ND performs slightly better than AMD.  In fact, for a 
fine mesh using $n = 20$, ND reduced fill-in by a factor of 2.3.  In contrast, 
for the same mesh, AMD reduces fill-in by a factor of 1.9 and RCM reduces 
fill-in by a factor of 1.5.

\begin{table}[!htb]
\centering
\begin{tabular}{|c||c|c|c|}
\hline
n & Total Entries & Nonzero Entries & LU \\
\hline \hline
5 & 22500 & 1620 & 4680 \\
\hline
10 & $3.60 \times 10^5$ & 6840  & 36810 \\
\hline
15 & $1.82 \times 10^6$ & 15660 & $1.23 \times 10^5$ \\
\hline
20 & $5.76 \times 10^6$ & 28080 & $2.91 \times 10^5$ \\
\hline
\end{tabular} \\[.5cm]
\begin{tabular}{|c||c|c|c|}
\hline
n & AMD & ND & RCM \\
\hline \hline
5 & 2970 & 3258 & 3600 \\
\hline
10 & 20427 & 20508 & 26550 \\
\hline
15 & 69567 & 61647 & 87804 \\
\hline
20 & $1.53 \times 10^5$ & $1.33 \times 10^5$ & $2.04 \times 10^5$ \\
\hline
\end{tabular}
\caption{(Above) Total entries, number of nonzero entries before finding the 
$LU$ decomposition, and the number of nonzero entries in the $LU$ 
decomposition.
(Below) Number of nonzero entries of the $LU$ decomposition after applying AMD, 
ND, and RCM pivoting.\label{tab:Exp1_a}}
\end{table}

\begin{table}[!htb]
\centering
\begin{tabular}{|c||c|c|c|c|c|}
\hline
n & Nonzero Entries & LU & AMD & ND & RCM \\
\hline \hline
5 & 7.2 \% & 20.8 \% & 13.2 \% & 14.5 \% & 16 \% \\
\hline
10 & 1.9 \% & 10.2 \% & 5.7 \% & 5.7 \% & 7.4 \% \\
\hline
15 & 0.86 \% & 6.8 \% & 3.8 \% & 3.4 \% & 4.8 \% \\
\hline
20 & 0.49 \% & 5.1 \% & 2.7 \% & 2.3 \% & 3.5 \% \\
\hline
\end{tabular}
\caption{Percentage of nonzero entries before finding the $LU$ decomposition, 
in the $LU$ decomposition, and in the $LU$ decomposition after applying AMD, 
ND, and RCM pivoting. \label{tab:Exp1_b}}
\end{table}

\subsection{Numerical experiment 2} \label{sec:NumericalExp2}
\begin{figure}[!htb]
\includegraphics[width = .5\textwidth]{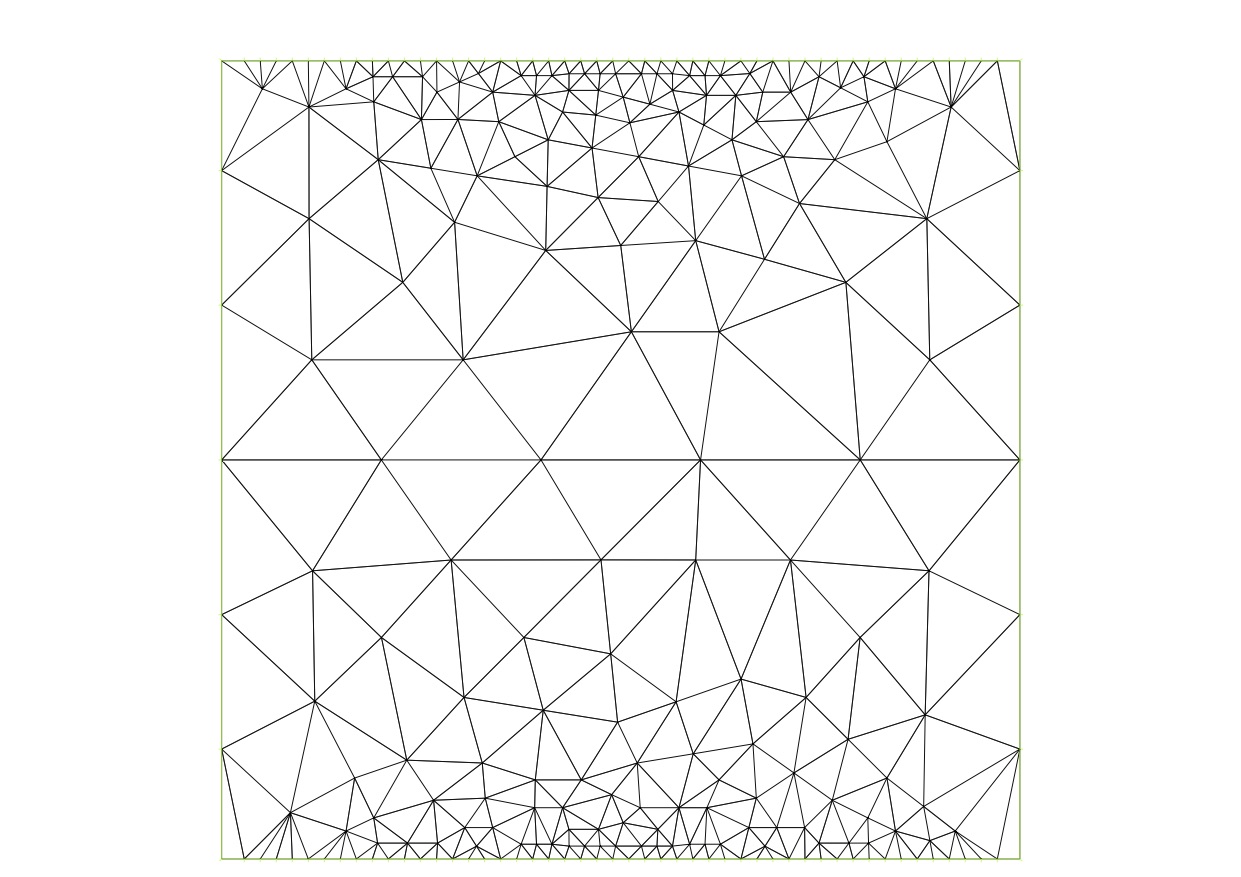}
\includegraphics[width = .5\textwidth]{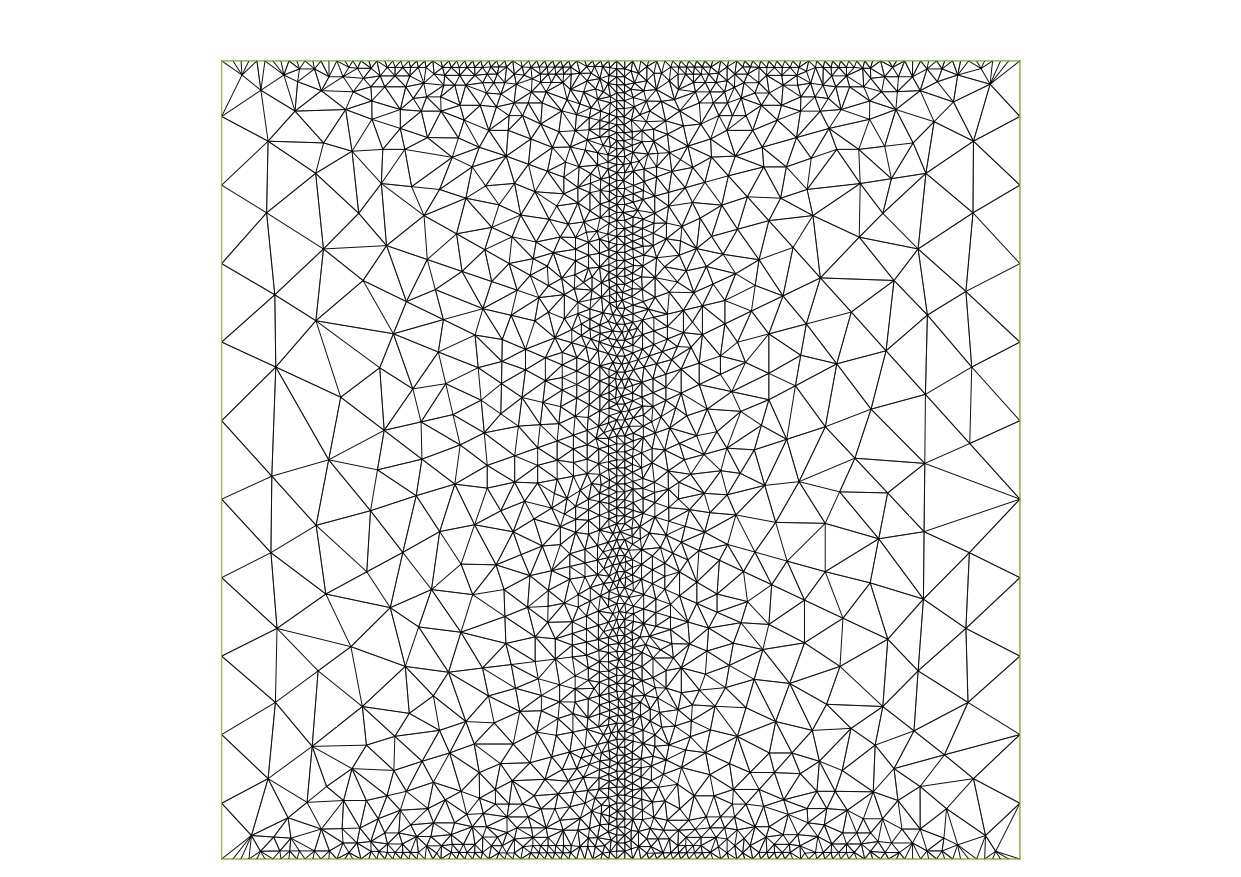}
\caption{\label{fig:Exp2Mesh} Mesh used in experiment 2 with $n = 5$ (left) 
and $n = 10$ (right).}
\end{figure}

Similar to the experiments in Section \ref{sec:NumericalExp1}, in this section 
we will let $D_2 = \emptyset$ and use the IP-DG function space $V^1_h$.  The 
main difference between the experiments in this section and Section 
\ref{sec:NumericalExp1} is the use of a less structured triangular mesh 
$\mathcal{T}_h$.  In this section we use a mesh $\mathcal{T}_h$ defined using 
$n = 5, 10, 15, 20$ intervals along the vertical sides of $D$ and $10n$ 
intervals along the horizontal sides of $D$.  Figure \ref{fig:Exp2Mesh} shows 
examples of meshes used in this section.  The structure of the mesh will 
determine the structure of the matrix $A$ in \eqref{eq:InnerProd}.  Thus, we 
expect to see different results in this experiment. Similar to Section 
\ref{sec:NumericalExp1} the IP-DG function space $V^1_h$ was used in this set 
of experiments.

\begin{figure}[!htb]
\includegraphics[width = .5\textwidth]{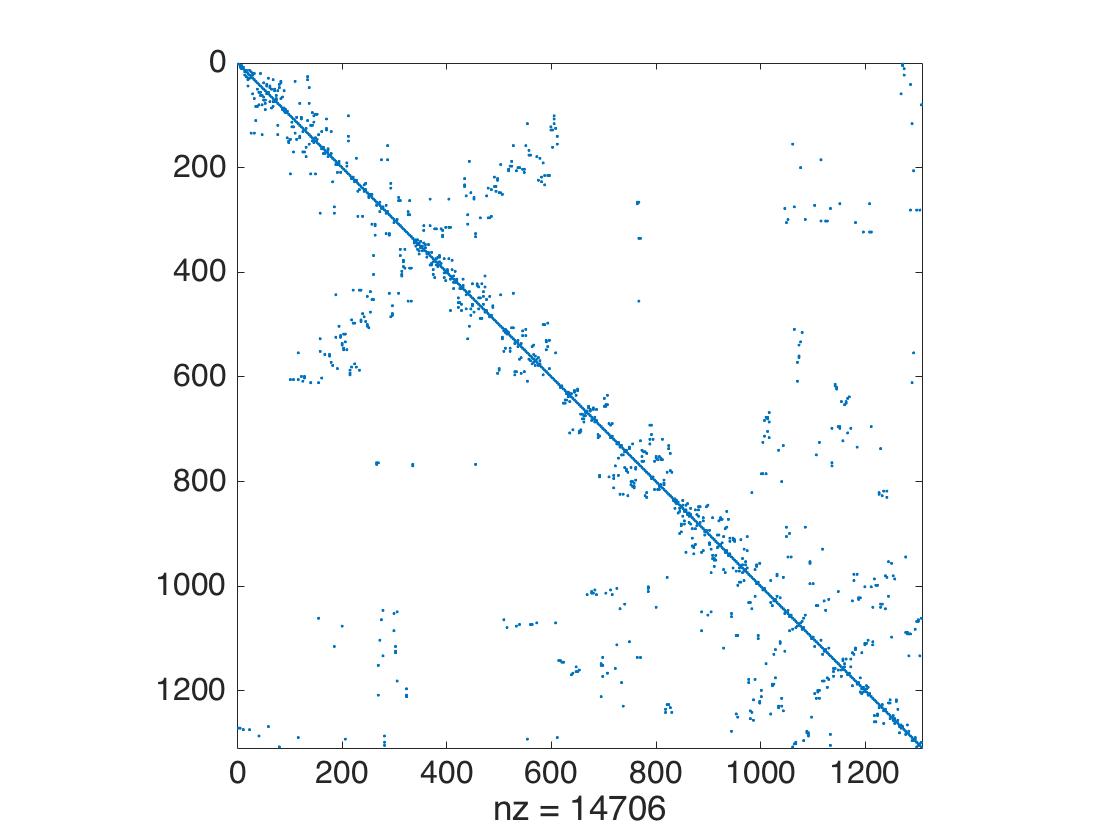}
\includegraphics[width = .5\textwidth]{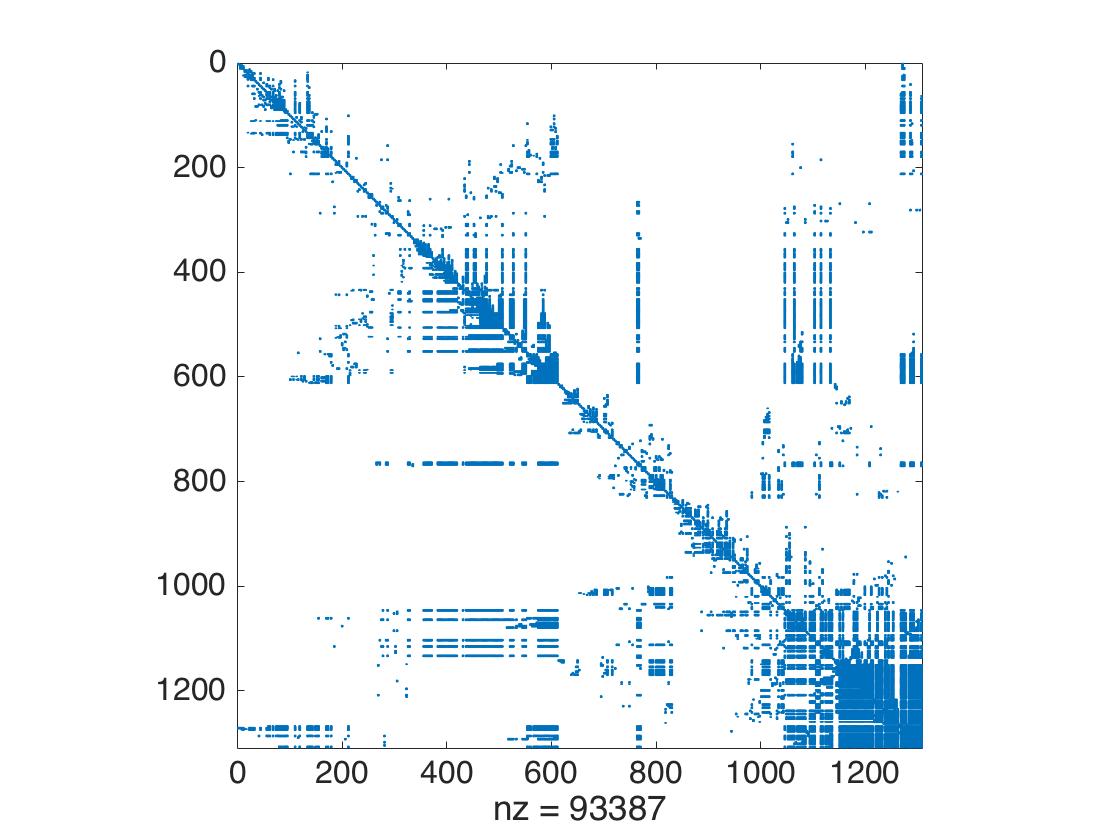}
\caption{\label{fig:Exp2NoPivot5} (Left) sparsity structure of the global 
matrix A produced by the IP-DG method with $n = 5$. (Right) sparsity structure 
of the combined LU decomposition of A.}
\end{figure}

\begin{figure}[!htb]
\includegraphics[width = .5\textwidth]{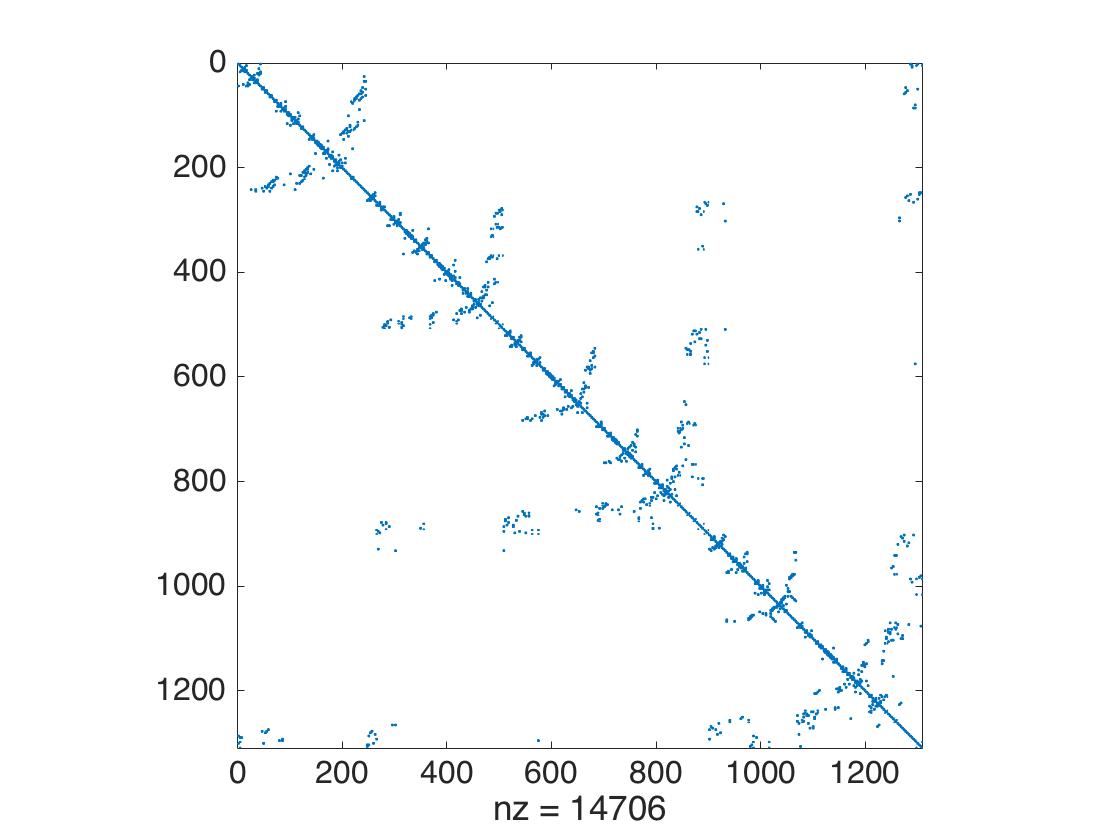}
\includegraphics[width = .5\textwidth]{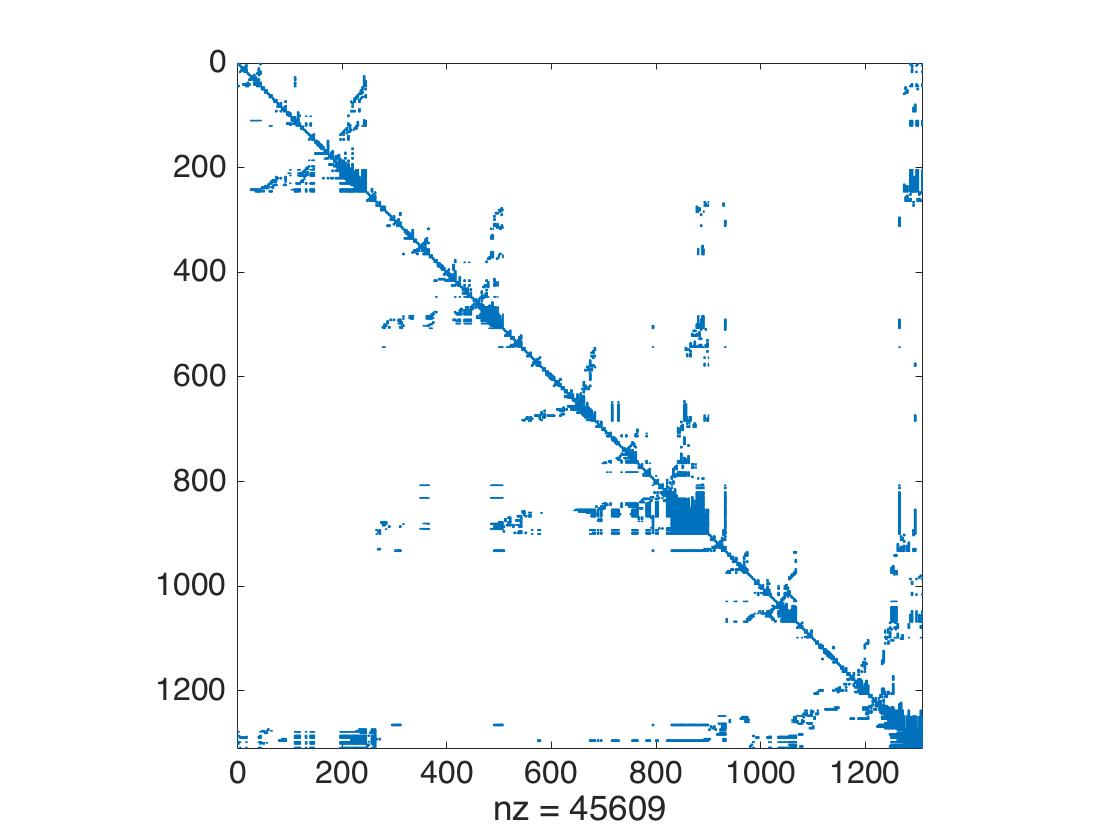}
\caption{\label{fig:Exp2AMD5} (Left) sparsity structure of the global matrix 
after AMD pivoting. (Right) sparsity structure of the combined LU 
decomposition the global matrix after AMD pivoting.}
\end{figure}

\begin{figure}[!htb]
\includegraphics[width = .5\textwidth]{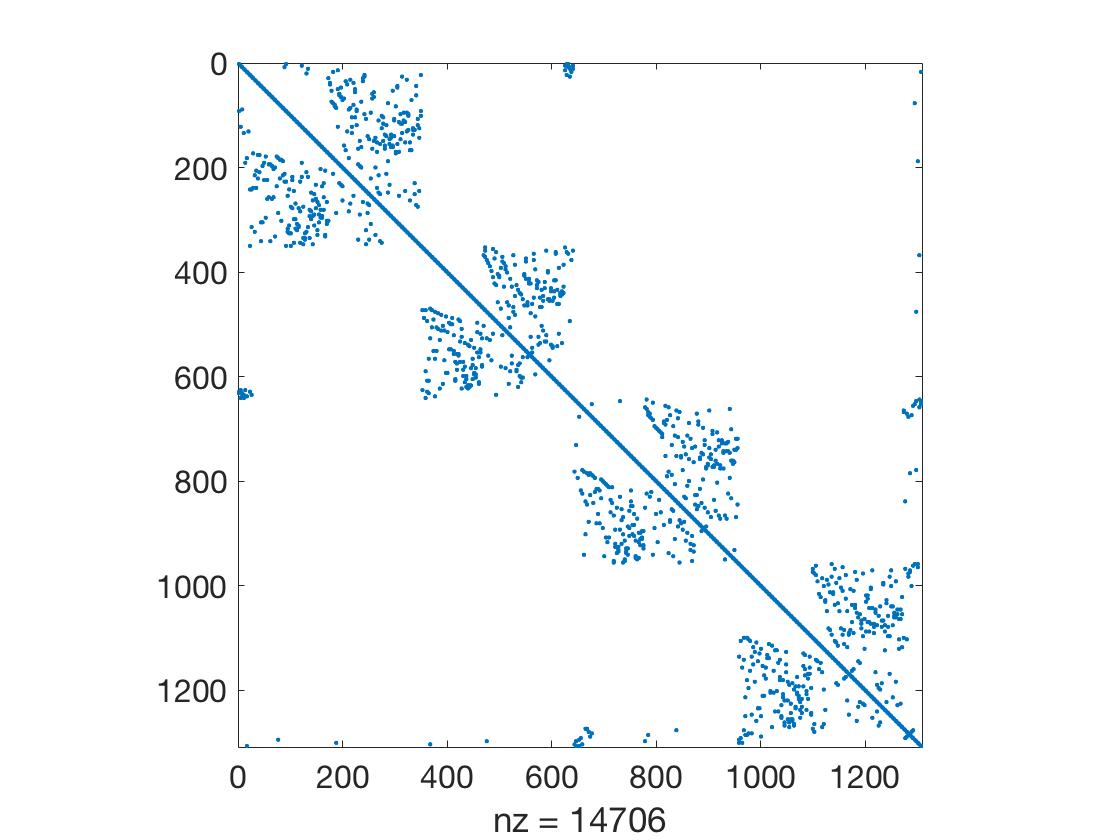}
\includegraphics[width = .5\textwidth]{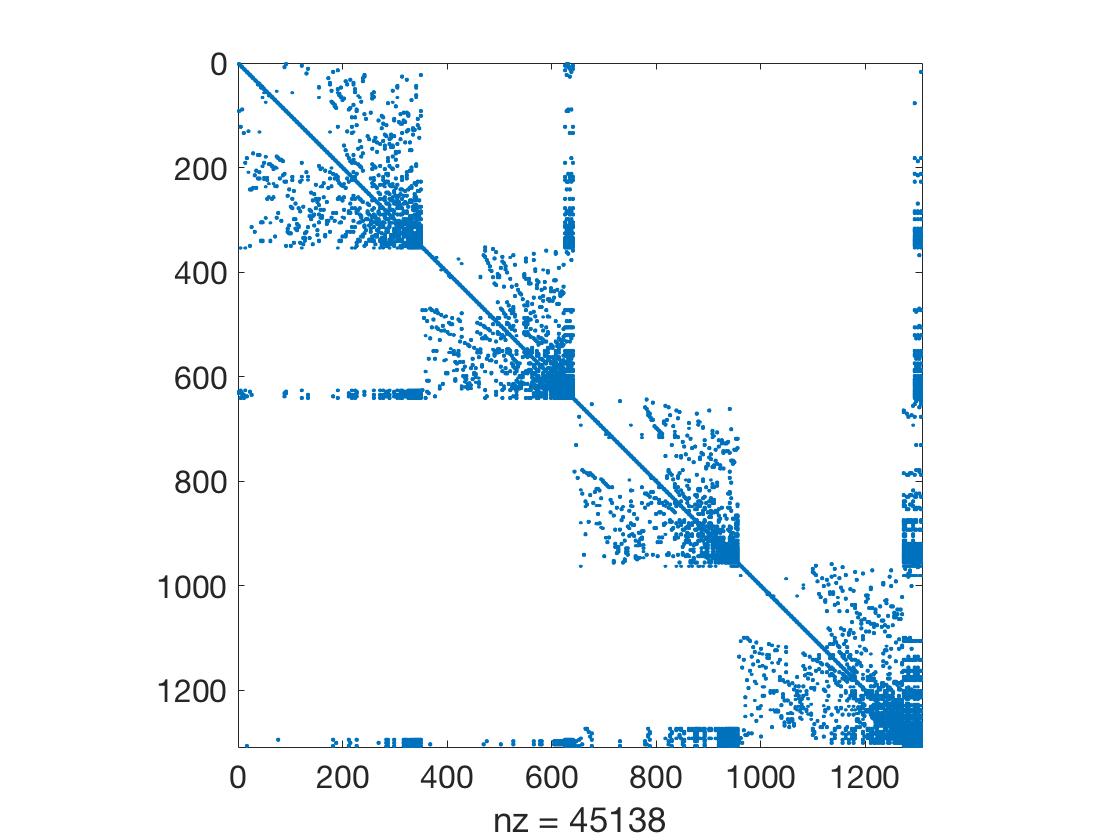}
\caption{\label{fig:Exp2ND5} (Left) sparsity structure of the global matrix 
after ND pivoting. (Right) sparsity structure of the combined LU decomposition 
the global matrix after ND pivoting.}
\end{figure}

\begin{figure}[!htb]
\includegraphics[width = .5\textwidth]{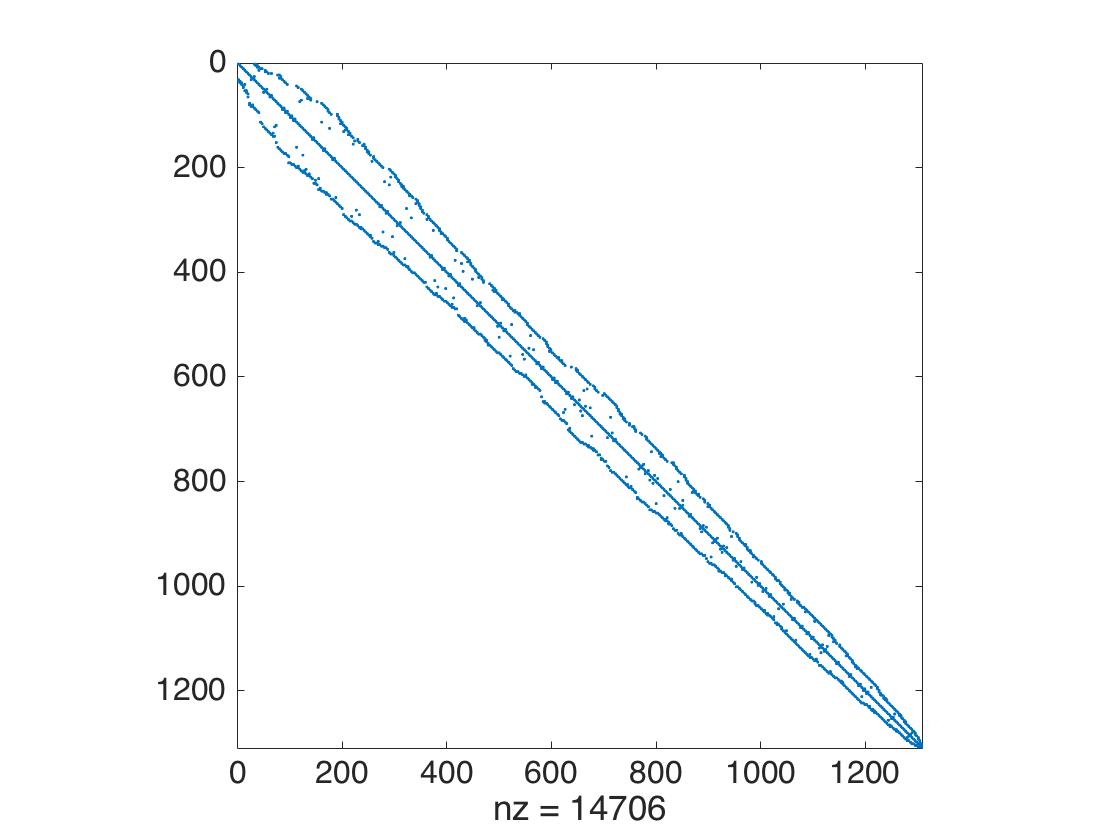}
\includegraphics[width = .5\textwidth]{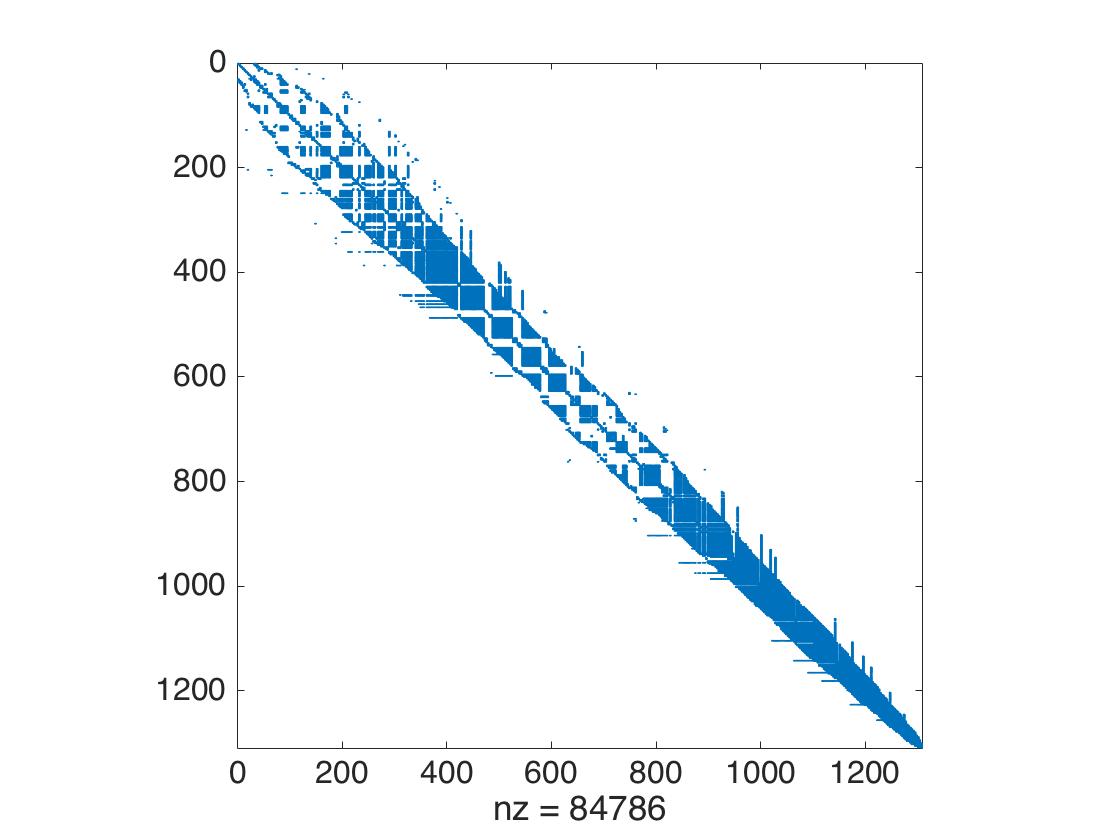}
\caption{\label{fig:Exp2RCM5} (Left) sparsity structure of the global matrix 
after RCM pivoting. (Right) sparsity structure of the combined LU 
decomposition the global matrix after RCM pivoting.}
\end{figure}

Similar to Section \ref{sec:NumericalExp1}, we plotted the sparsity structure 
of the matrix $A$ and its combined $LU$ factorization for $n = 5$.  We also 
plotted the matrix $PAP^T$ after AMD, ND, and RCM pivoting along with the $LU$ 
factorization of each of these pivoted matrices.  These plots are given in 
Figure \ref{fig:Exp2NoPivot5}--\ref{fig:Exp2RCM5}.  Also, similar to Section 
\ref{sec:NumericalExp1} the data for the number of nonzero entries for the 
matrices in this set of experiments are summarized in Tables \ref{tab:Exp2_a} 
and \ref{tab:Exp2_b}.  From these tables we see that all three pivoting 
strategies reduce the fill-in associated to $LU$ factorization in this set of  
experiments. Also, as the mesh parameter $n$ increases, the pivoting 
strategies perform better.  In this experiment, RCM reduces fill-in the least 
and ND reduces fill-in the most.  For instance, for $n = 20$, RCM reduces 
fill-in by a factor of $1.99$ and ND reduces fill-in by a factor of $5.36$.

\begin{table}[!htb]
\centering
\begin{tabular}{|c||c|c|c|}
\hline
n & Total Entries & Nonzero Entries & LU \\
\hline \hline
5 & $1.71 \times 10^6$ & 14706 & 93387 \\
\hline
10 & $9.91 \times 10^7$ & $1.17 \times 10^5$ & $4.52 \times 10^6$ \\
\hline
15 & $1.24 \times 10^8$ & $1.31 \times 10^5$ & $3.24 \times 10^6$ \\
\hline
20 & $1.46 \times 10^9$ & $4.55 \times 10^5$ & $5.88 \times 10^7$ \\
\hline
\end{tabular}\\[.5cm]
\begin{tabular}{|c||c|c|c|}
\hline
n & AMD & ND & RCM \\
\hline \hline
5 & 45609 & 45138 & 84786 \\
\hline
10 & $1.41 \times 10^6$ & $9.71 \times 10^5$ & $3.33 \times 10^6$ \\
\hline
15 & $1.43 \times 10^6$ & $1.24 \times 10^6$ & $2.79 \times 10^6$ \\
\hline
20 & $1.97 \times 10^7$ & $1.10 \times 10^7$ & $2.95 \times 10^7$ \\
\hline
\end{tabular}
\caption{(Above) Total entries, number of nonzero entries before finding the 
$LU$ decomposition, and the number of nonzero entries in the $LU$ 
decomposition. (Below) Number of nonzero entries of the $LU$ decomposition 
after applying AMD, ND, and RCM pivoting.\label{tab:Exp2_a}}
\end{table}

\begin{table}[!htb]
\centering
\begin{tabular}{|c||c|c|c|c|c|}
\hline
n & Nonzero Entries & LU & AMD & ND & RCM \\
\hline \hline
5 & 0.86 \% & 5.5 \% & 2.7 \% & 2.6 \% & 4.96 \% \\
\hline
10 & 0.12 \% & 4.6 \% & 1.4 \% & 0.98 \% & 3.4 \% \\
\hline
15 & 0.11 \% & 2.6 \% & 1.2 \% & 1.0 \% & 2.2 \% \\
\hline
20 & 0.031 \% & 4.02 \% & 1.3 \% & 0.75 \% & 2.02 \% \\
\hline
\end{tabular}
\caption{Percentage of nonzero entries before finding the $LU$ decomposition, 
in the $LU$ decomposition, and in the $LU$ decomposition after applying AMD, 
ND, and RCM pivoting.\label{tab:Exp2_b}}
\end{table}

\clearpage

\subsection{Numerical experiment 3} \label{sec:NumericalExp3}
\begin{figure}[!htb]
\includegraphics[width = .5\textwidth]{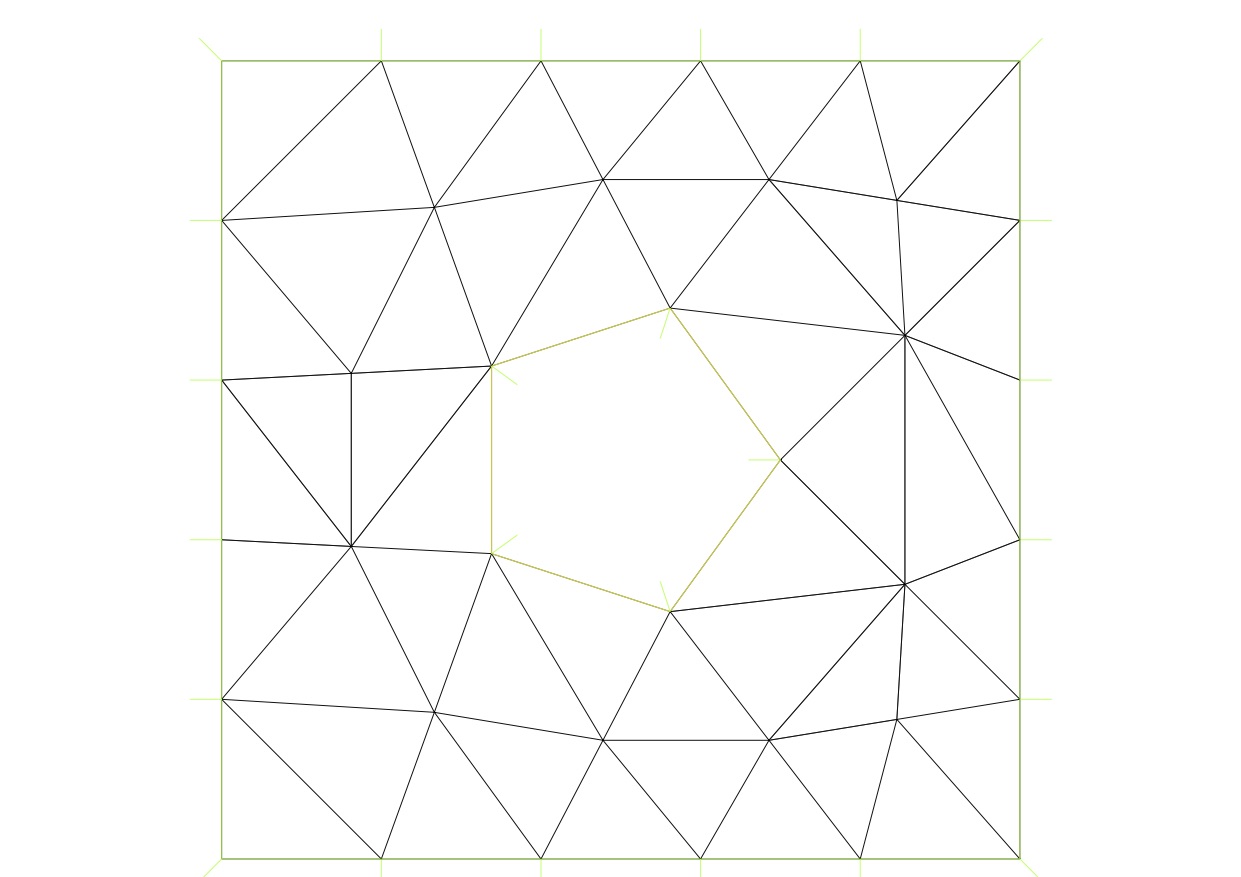}
\includegraphics[width = .5\textwidth]{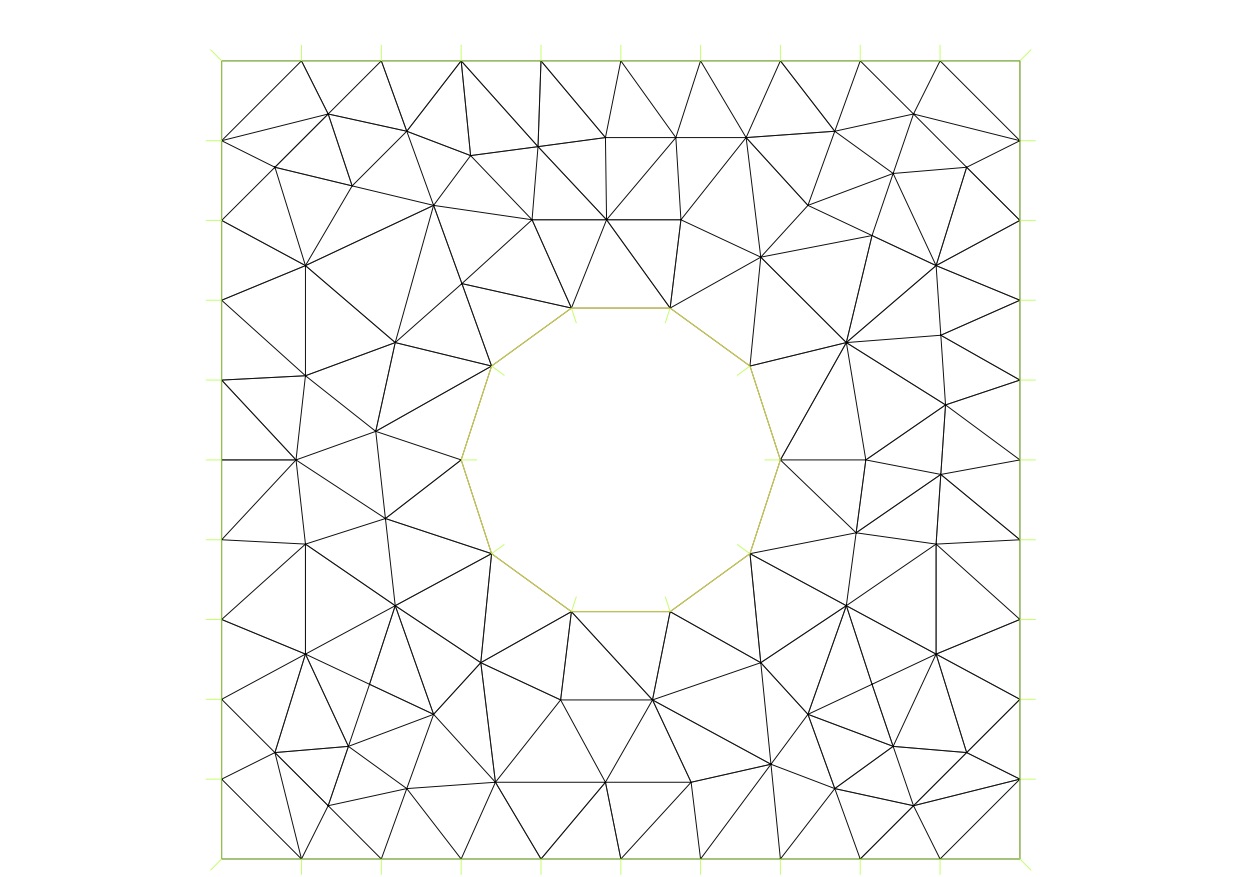}
\caption{\label{fig:Exp3Mesh} Mesh used in experiment 3 with $n = 5$ (left) 
and $n = 10$ (right).}
\end{figure}

In the third set of numerical experiments presented in this paper, our goal 
was to observe the performance of sparsity preserving pivoting on a domain 
with a scattering boundary.  With this goal in mind, a circular scattering 
domain
\begin{align}
	D_2 = \left\{ \, (x,y) \, \big| \, x^2 + y^2 \leq 1/5 \, \right\},
	\label{eq:ScatteringDomain}
\end{align}
was included.  In this experiment we used a quasi-uniform triangular mesh 
$\mathcal{T}_h$ generated by using $n = 5, 10, 15, 20$ intervals along each 
side of $\Gamma_1 = \pa D_1$ and along $\Gamma_2 = \pa D_2$.  Examples of the 
mesh generated by $n = 5, 10$ are given in Figure \ref{fig:Exp3Mesh}.
Similar to Sections \ref{sec:NumericalExp1} and \ref{sec:NumericalExp2} the 
IP-DG function space $V^1_h$ was used in this set of experiments.

\begin{figure}[!htb]
\includegraphics[width = .5\textwidth]{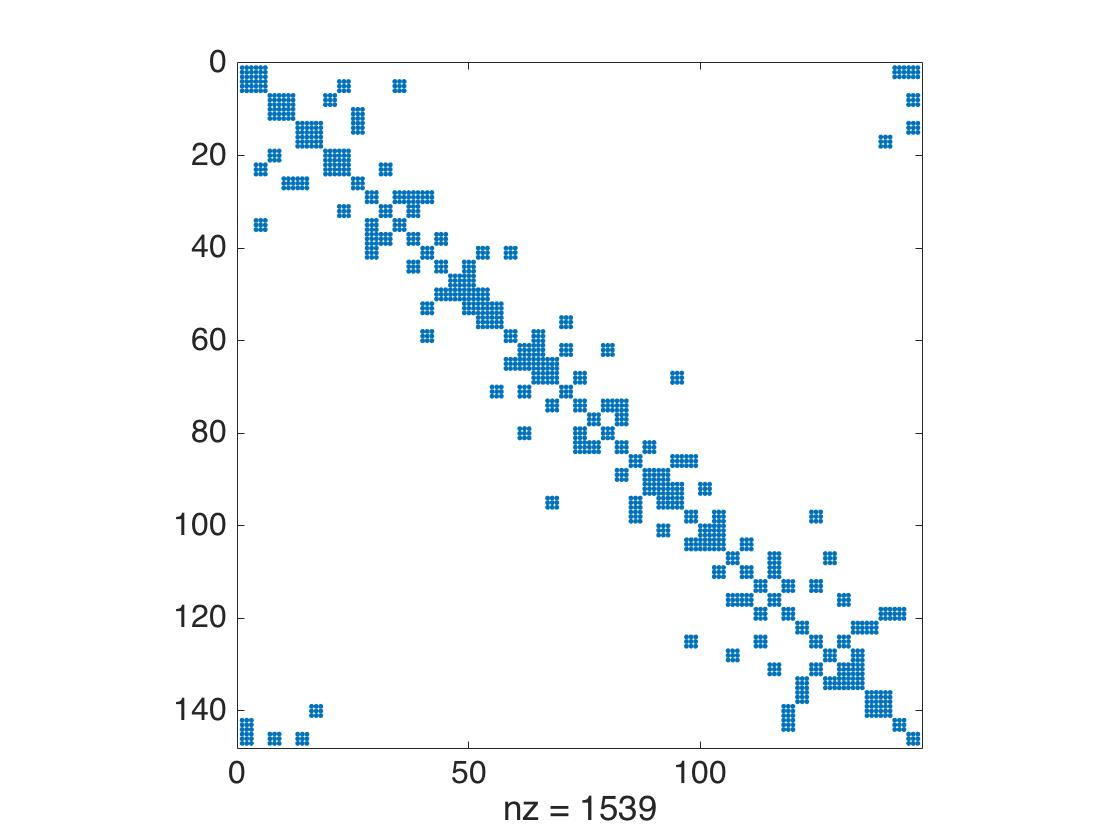}
\includegraphics[width = .5\textwidth]{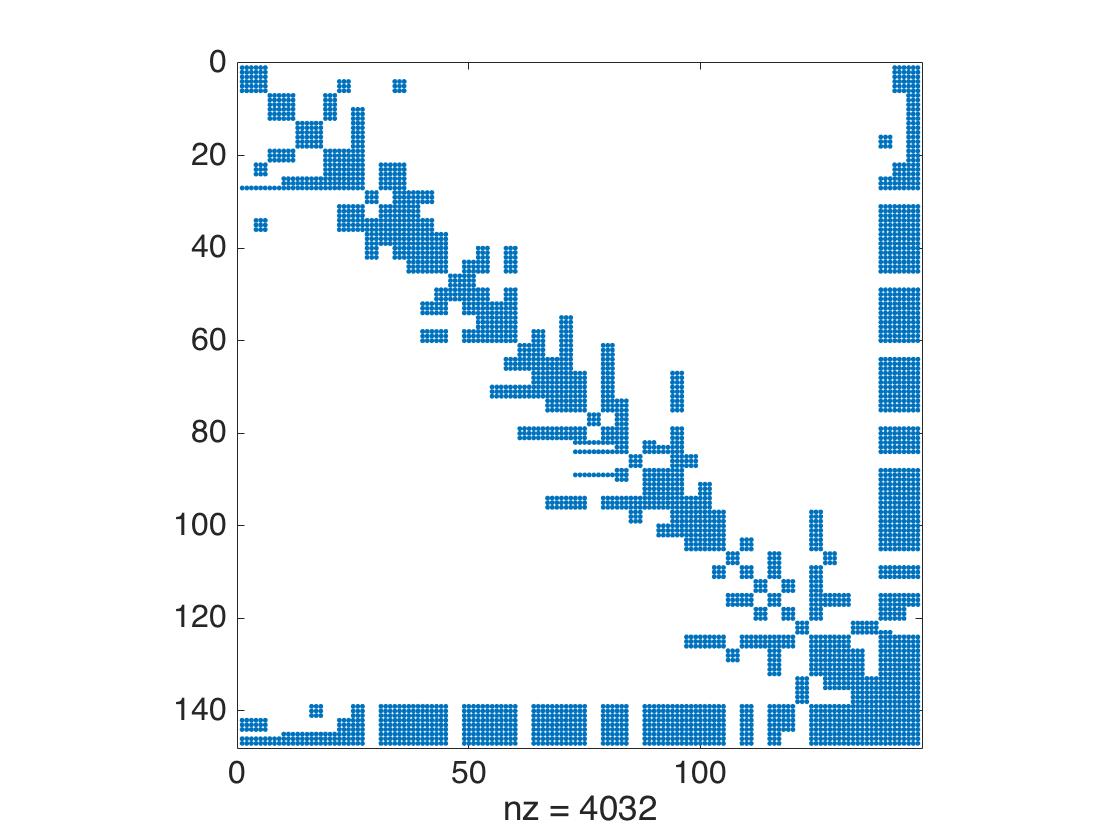}
\caption{\label{fig:Exp3NoPivot5} (Left) sparsity structure of the global 
matrix A produced by the IP-DG method with $n = 5$. (Right) sparsity structure 
of the combined LU decomposition of A.}
\end{figure}

\begin{figure}[!htb]
\includegraphics[width = .5\textwidth]{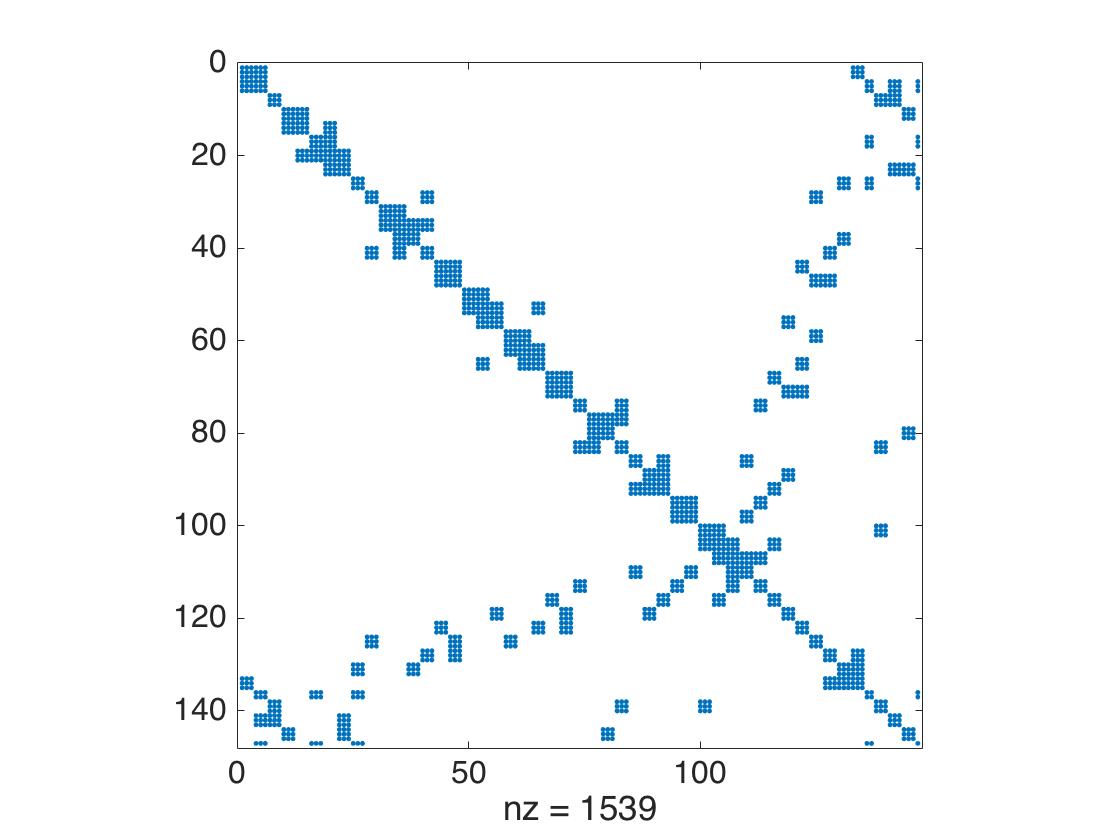}
\includegraphics[width = .5\textwidth]{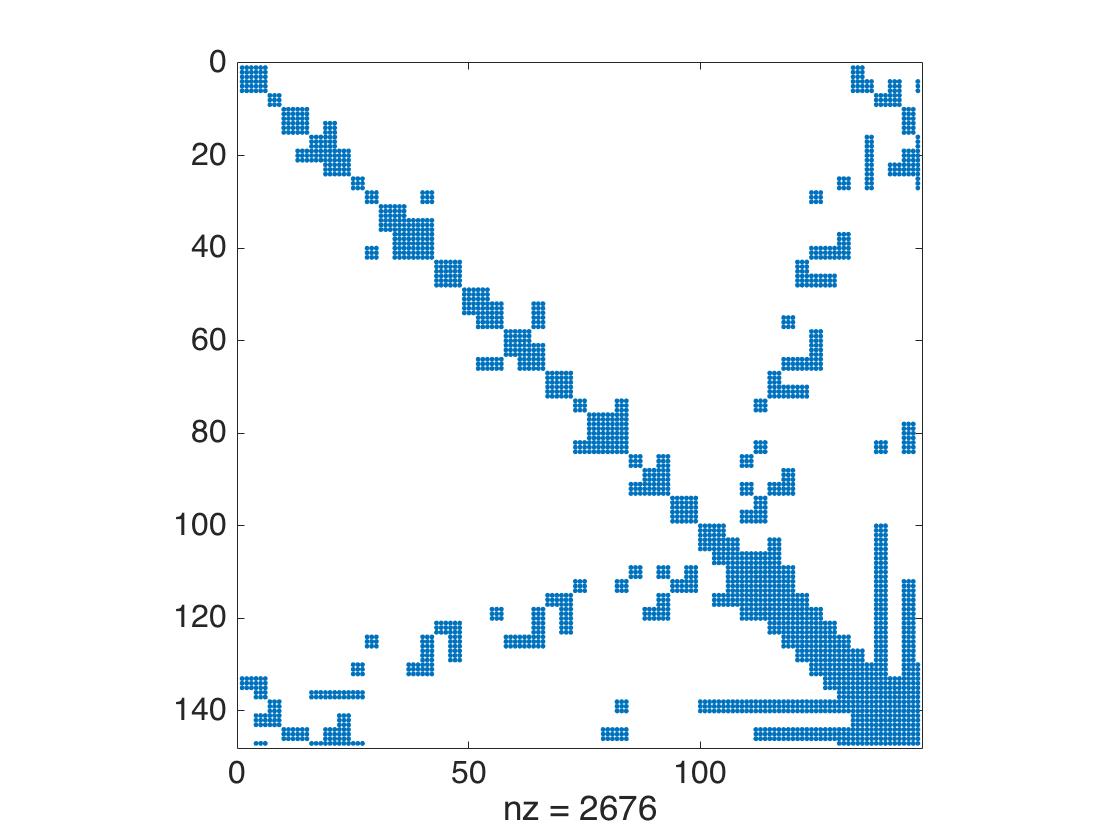}
\caption{\label{fig:Exp3AMD5} (Left) sparsity structure of the global matrix 
after AMD pivoting. (Right) sparsity structure of the combined LU 
decomposition the global matrix after AMD pivoting.}
\end{figure}

\begin{figure}[!htb]
\includegraphics[width = .5\textwidth]{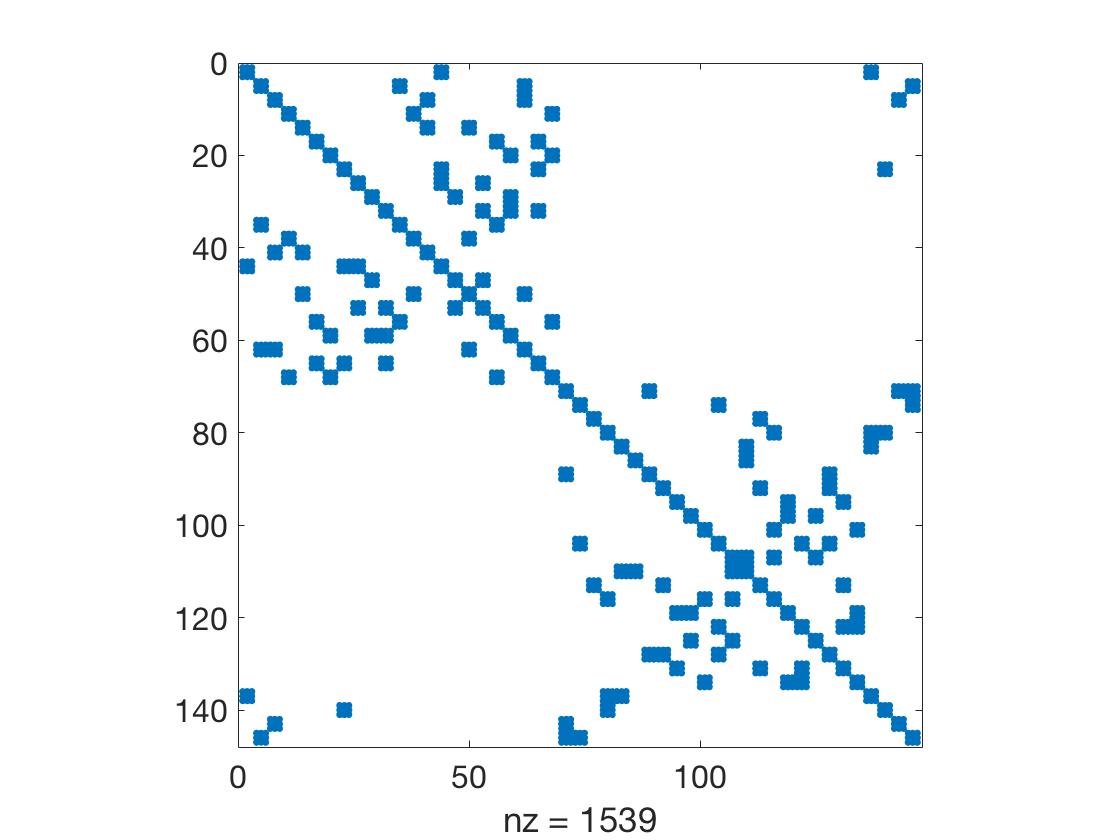}
\includegraphics[width = .5\textwidth]{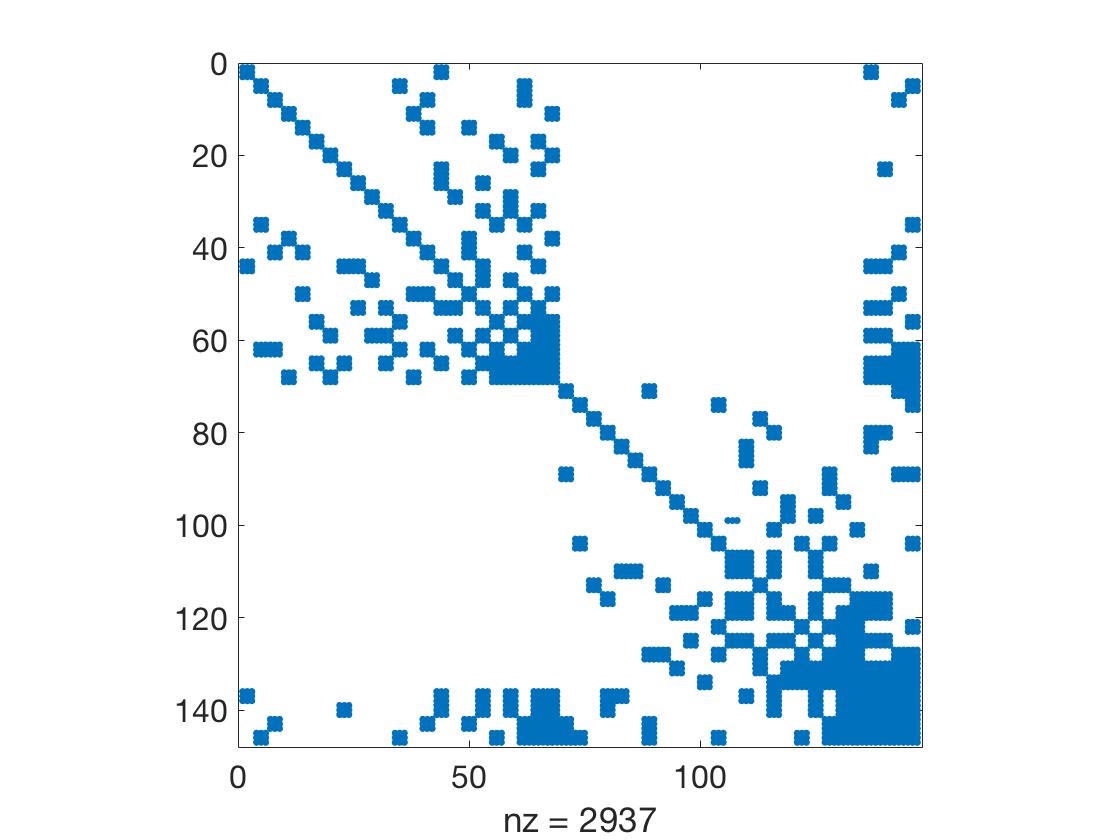}
\caption{\label{fig:Exp3ND5} (Left) sparsity structure of the global matrix 
after ND pivoting. (Right) sparsity structure of the combined LU decomposition 
the global matrix after ND pivoting.}
\end{figure}

\begin{figure}[!htb]
\includegraphics[width = .5\textwidth]{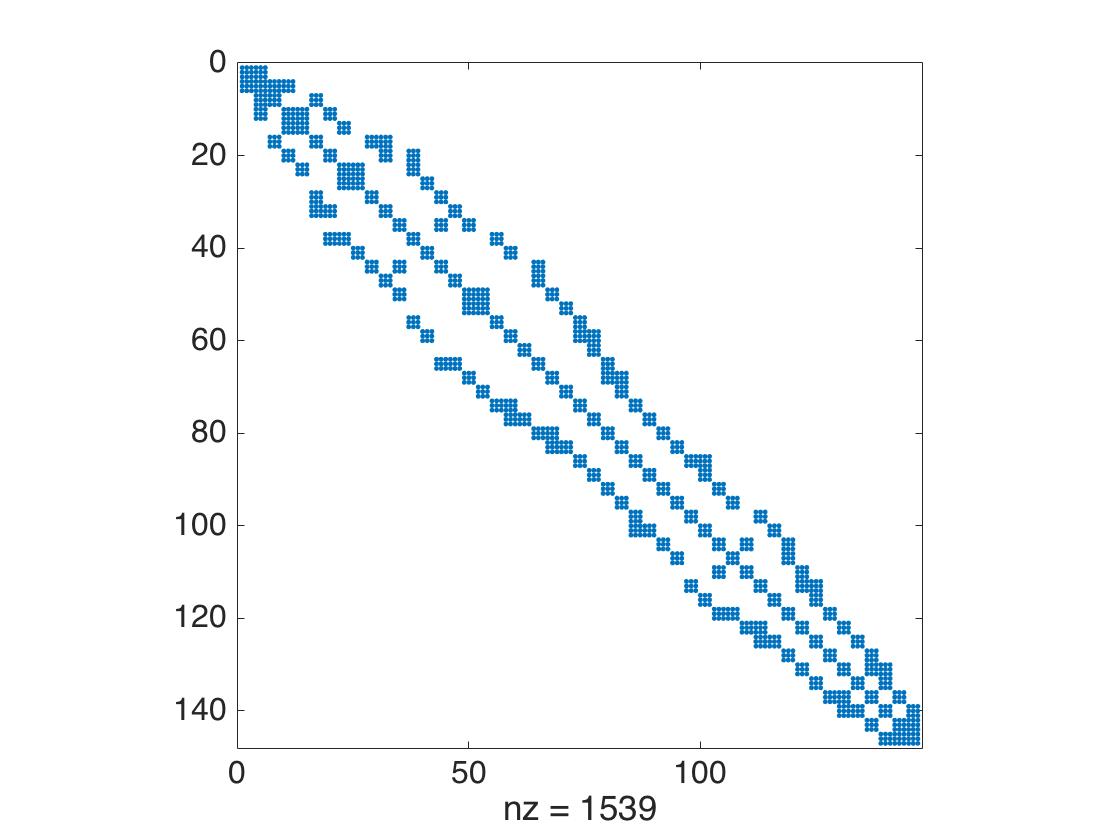}
\includegraphics[width = .5\textwidth]{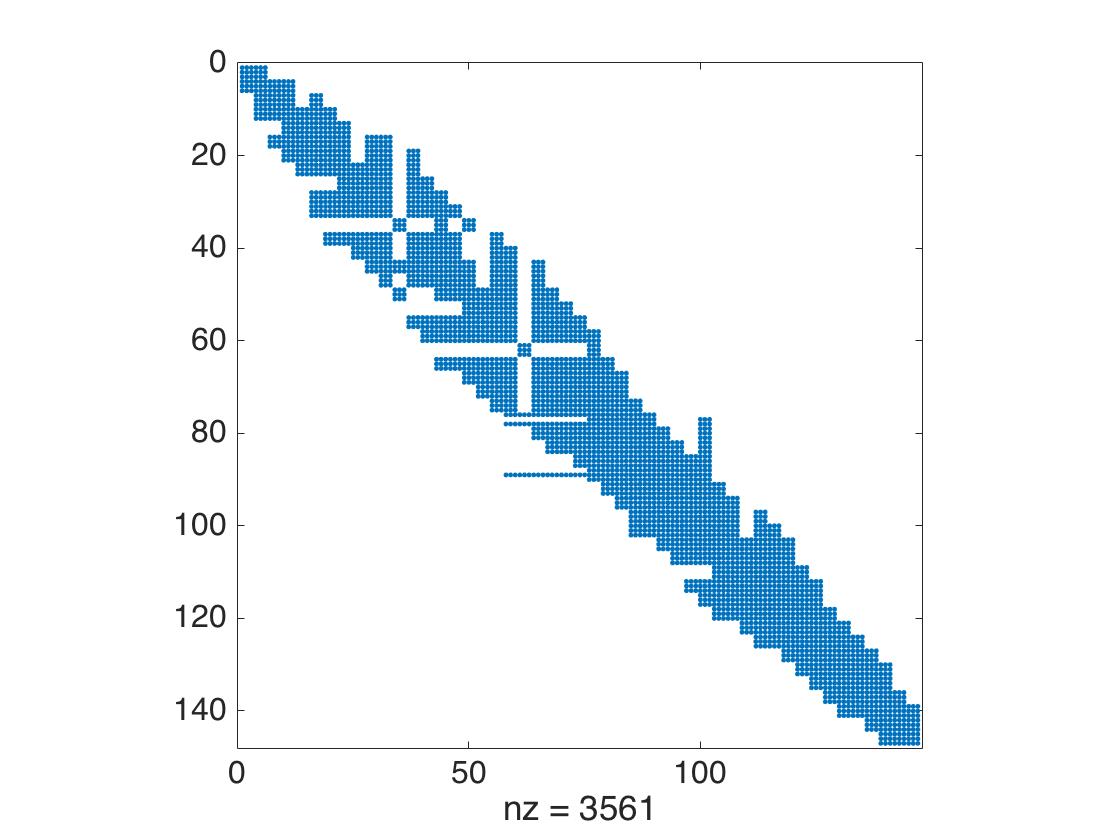}
\caption{\label{fig:Exp3RCM5} (Left) sparsity structure of the global matrix 
after RCM pivoting. (Right) sparsity structure of the combined LU 
decomposition the global matrix after RCM pivoting.}
\end{figure}

Figure \ref{fig:Exp3NoPivot5} shows the sparsity structure of the matrix $A$ 
defined by \eqref{eq:LinSystem} along with its combined $LU$ factorization.  
Similarly, Figures \ref{fig:Exp3AMD5}--\ref{fig:Exp3RCM5} give the sparsity 
structure of $PAP^T$ after AMD, ND, and RCM pivoting along with the combined 
$LU$ factorization of each matrix for $n = 5$.  The number of nonzero entries 
for $A$ and the combined $LU$ factorization of $A$ and $PAP^T$ after AMD, ND, 
and RCM are summarized in Table \ref{tab:Exp3_a}.  Similarly, the percent of 
non-zero entries for the combined $LU$ factorization of these matrices is 
given in Table \ref{tab:Exp3_b}.


\begin{table}[!htb]
\centering
\begin{tabular}{|c||c|c|c|}
\hline
n & Total Entries & Nonzero Entries & LU \\
\hline \hline
5 & 21609 & 1539  & 4032 \\
\hline
10 & $2.60 \times 10^5$ & 5670  & 24749 \\
\hline
15 & $1.25 \times 10^6$ & 12753 & $1.17 \times 10^5$ \\
\hline
20 & $3.99 \times 10^6$ & 23076 & $3.19 \times 10^5$ \\
\hline
\end{tabular} \\[.5cm]
\begin{tabular}{|c||c|c|c|}
\hline
n & AMD & ND & RCM \\
\hline \hline
5 & 2676 & 2937 & 3561 \\
\hline
10 & 17676 & 18429 & 21675 \\
\hline
15 & 49179 & 47826 & 77361 \\
\hline
20 & $1.29 \times 10^5$ & $1.26 \times 10^5$ & $1.92 \times 10^5$ \\
\hline
\end{tabular}
\caption{(Above) Total entries, number of nonzero entries before finding the 
$LU$ decomposition, and the number of nonzero entries in the $LU$ 
decomposition.
(Below) Number of nonzero entries of the $LU$ decomposition after applying AMD, 
ND, and RCM pivoting.\label{tab:Exp3_a}}
\end{table}

\begin{table}[!htb]
\centering
\begin{tabular}{|c||c|c|c|c|c|}
\hline
n & Nonzero Entries & LU & AMD & ND & RCM \\
\hline \hline
5 & 7.1 \% & 18.7 \% & 12.4 \% & 13.6 \% & 16.5 \% \\
\hline
10 & 2.2 \% & 9.5 \% & 6.8 \% & 7.1 \% & 8.3 \% \\
\hline
15 & 1.02 \% & 9.4 \% & 3.9 \% & 3.8 \% & 6.2 \% \\
\hline
20 & 0.58 \% & 7.99 \% & 3.2 \% & 3.1 \% & 4.8 \% \\
\hline
\end{tabular}
\caption{Percentage of nonzero entries before finding the $LU$ decomposition, 
in the $LU$ decomposition, and in the $LU$ decomposition after applying AMD, 
ND, and RCM pivoting.\label{tab:Exp3_b}}
\end{table}

The results of this experiment and summarized in Table \ref{tab:Exp3_b} are 
similar to the results using the uniform mesh in Section 
\ref{sec:NumericalExp1}.  For instance, Table \ref{tab:Exp3_b} indicates that 
AMD and ND reduce fill-in at roughly the same rate, with AMD performing 
slightly better for the coarse mesh $n = 5$ and ND performing slightly better 
for the fine mesh $n = 20$.  Similar to experiments in Section 
\ref{sec:NumericalExp1} and \ref{sec:NumericalExp2} all three pivoting 
strategies reduce fill-in associated to $LU$ factorization and as $n$ 
increases the rate at which pivoting reduces fill-in increases.  For the fine 
mesh given by $n = 20$, RCM reduces fill-in the least with a reduction factor 
of $1.66$ and ND reduces fill-in the most with a reduction factor of $2.58$.

\subsection{Numerical experiment 4} \label{sec:NumericalExp4}

\begin{figure}[!htb]
\includegraphics[width = .5\textwidth]{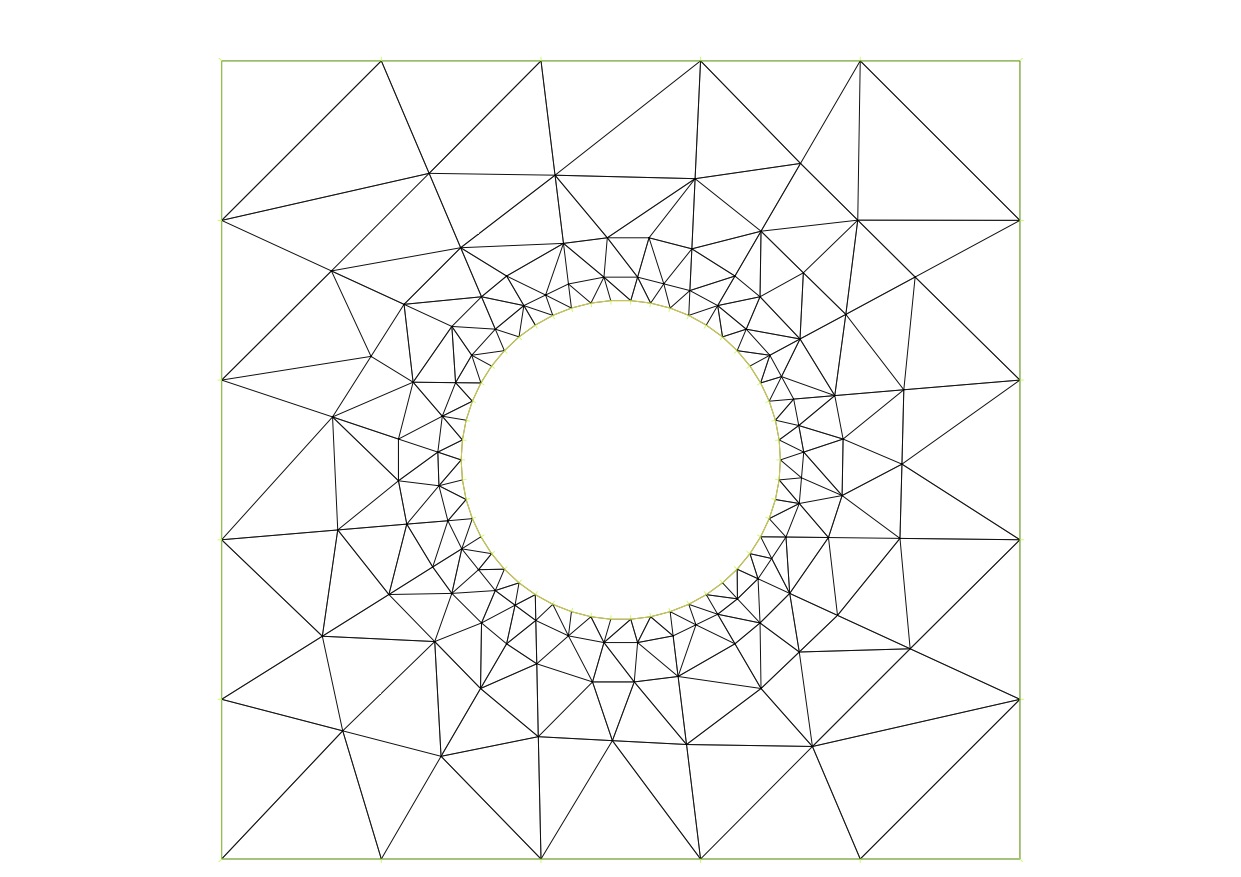}
\includegraphics[width = .5\textwidth]{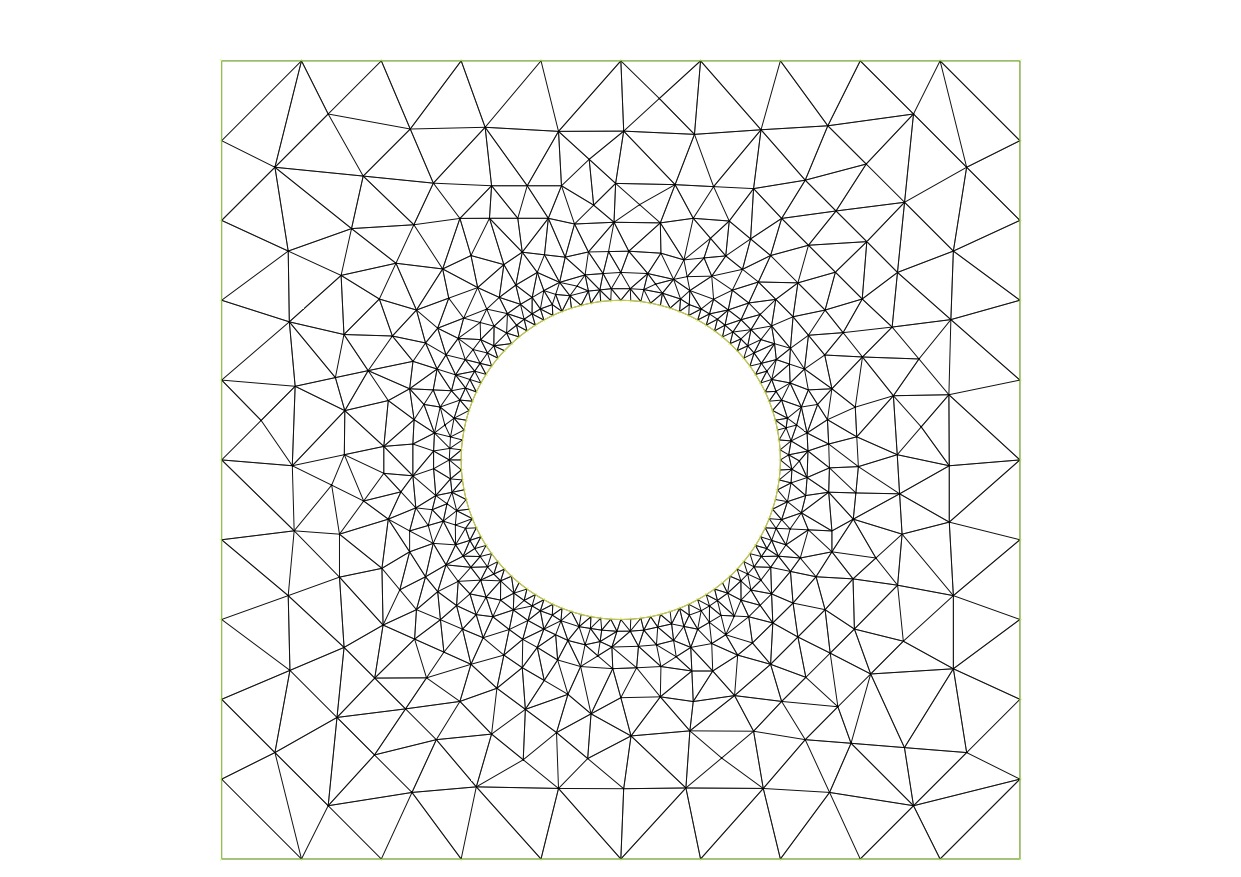}
\caption{\label{fig:Exp4Mesh} Mesh used in experiment 4 with $n = 5$ (left) 
and $n = 10$ (right).}
\end{figure}

Similar to Section \ref{sec:NumericalExp3} this section studies the IP-DG 
method with a circular scattering domain $D_2$ define by 
\eqref{eq:ScatteringDomain}.  In this section we chose to use a non-uniform 
triangular mesh $\mathcal{T}_h$. The triangular mesh $\mathcal{T}_h$ was 
generated using $n = 5, 10, 15, 20$ intervals along each side of $\Gamma_1 = 
\pa D_1$ and $10n$ intervals along $\Gamma_2 = \pa D_2$.  Thus, 
$\mathcal{T}_h$ is coarse along the boundary $\Gamma_1$ and fine along the 
boundary $\Gamma_2$.  Example meshes for $n = 5,10$ are given in Figure 
\ref{fig:Exp4Mesh}. Similar to Sections 
\ref{sec:NumericalExp1}--\ref{sec:NumericalExp3} the IP-DG function space 
$V^1_h$ was used in this set of experiments.

\begin{figure}[!htb]
\includegraphics[width = .5\textwidth]{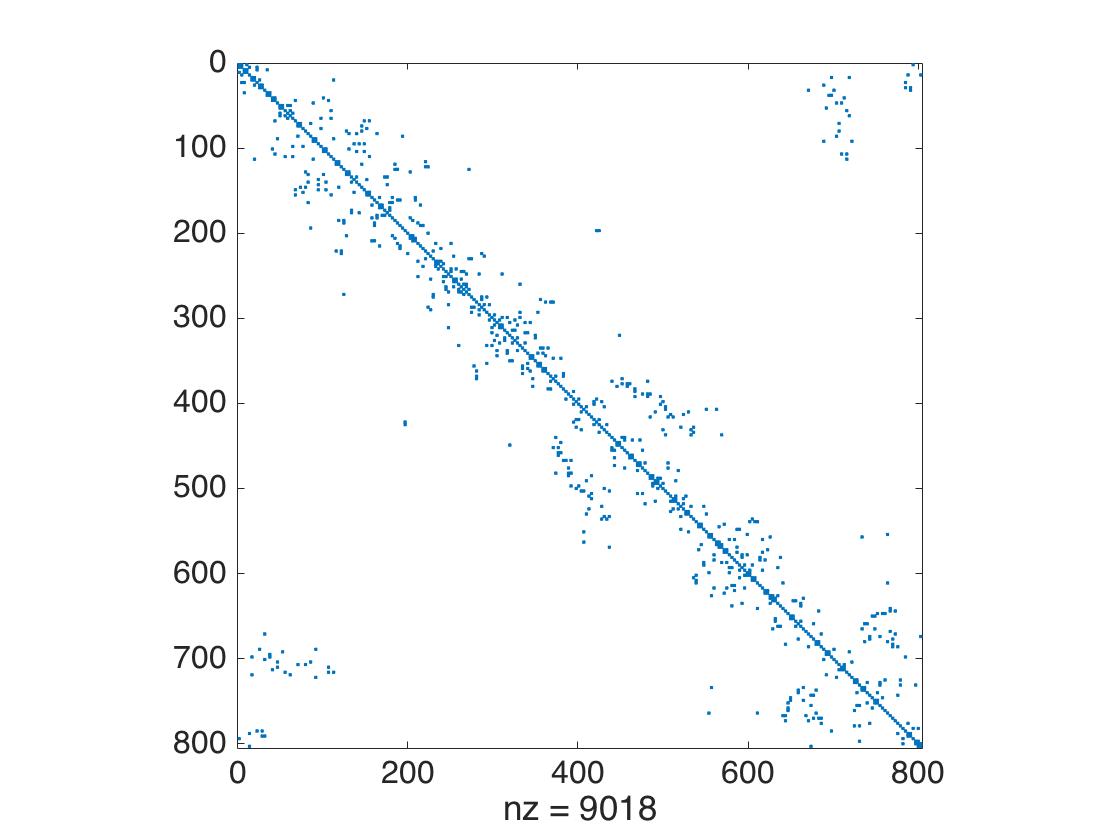}
\includegraphics[width = .5\textwidth]{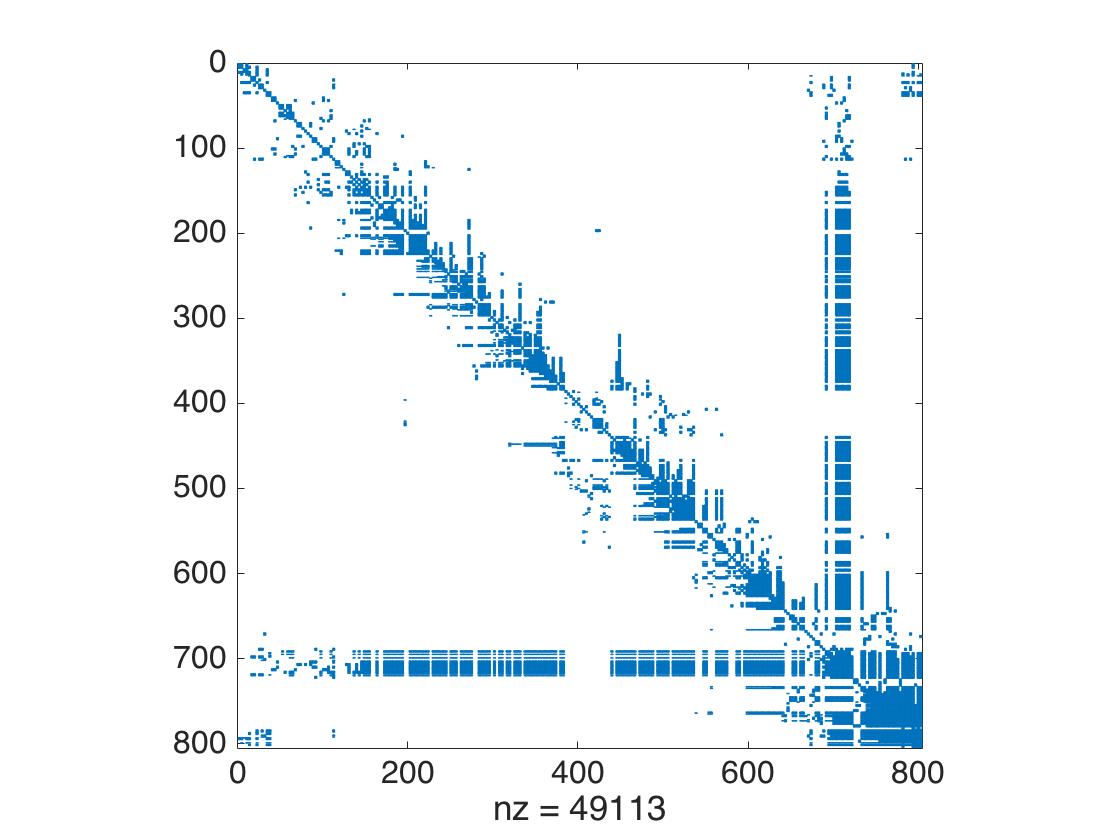}
\caption{\label{fig:Exp4NoPivot5} (Left) sparsity structure of the global 
matrix A produced by the IP-DG method with $n = 5$. (Right) sparsity structure 
of the combined LU decomposition of A.}
\end{figure}

\begin{figure}[!htb]
\includegraphics[width = .5\textwidth]{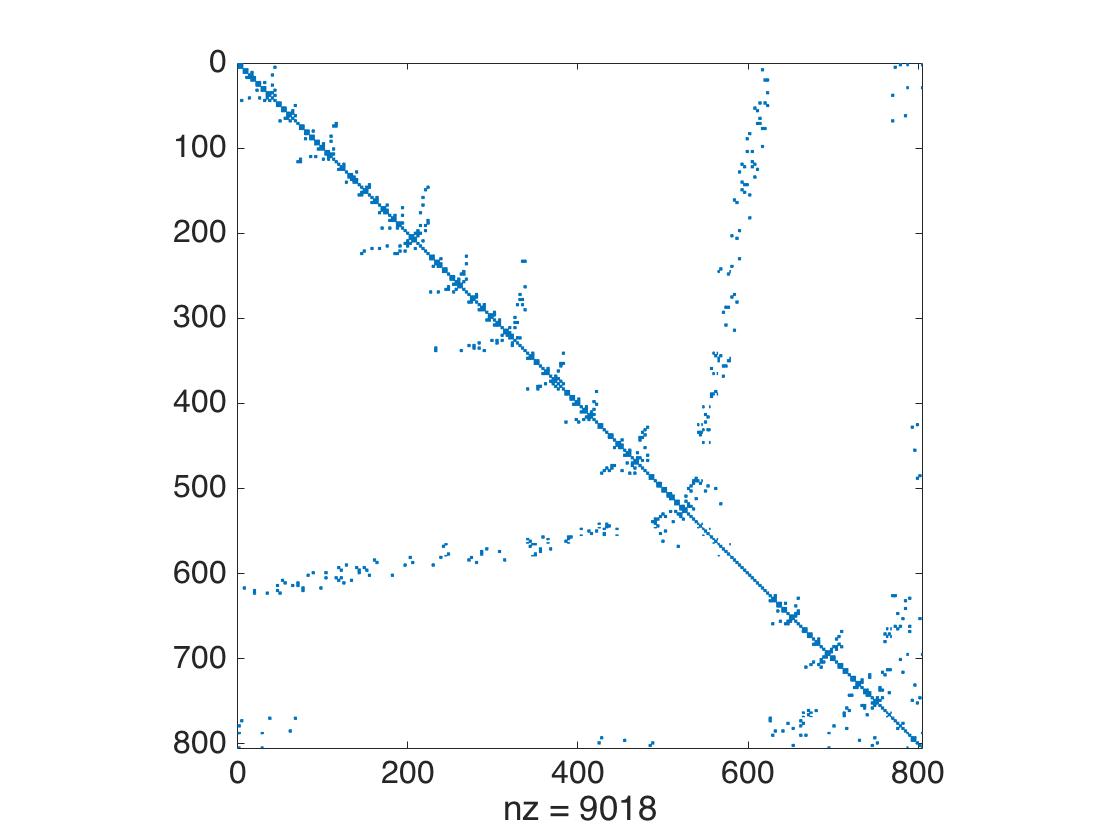}
\includegraphics[width = .5\textwidth]{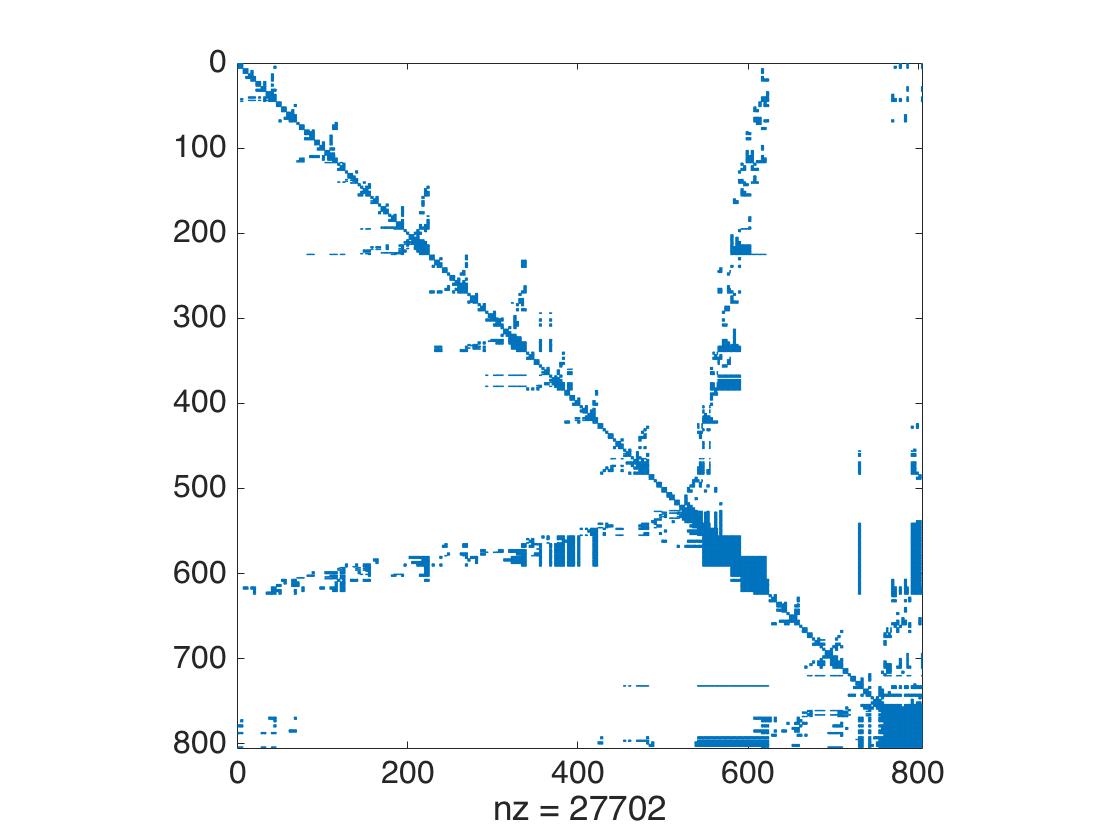}
\caption{\label{fig:Exp4AMD5} (Left) sparsity structure of the global matrix 
after AMD pivoting. (Right) sparsity structure of the combined LU 
decomposition the global matrix after AMD pivoting.}
\end{figure}

\begin{figure}[!htb]
\includegraphics[width = .5\textwidth]{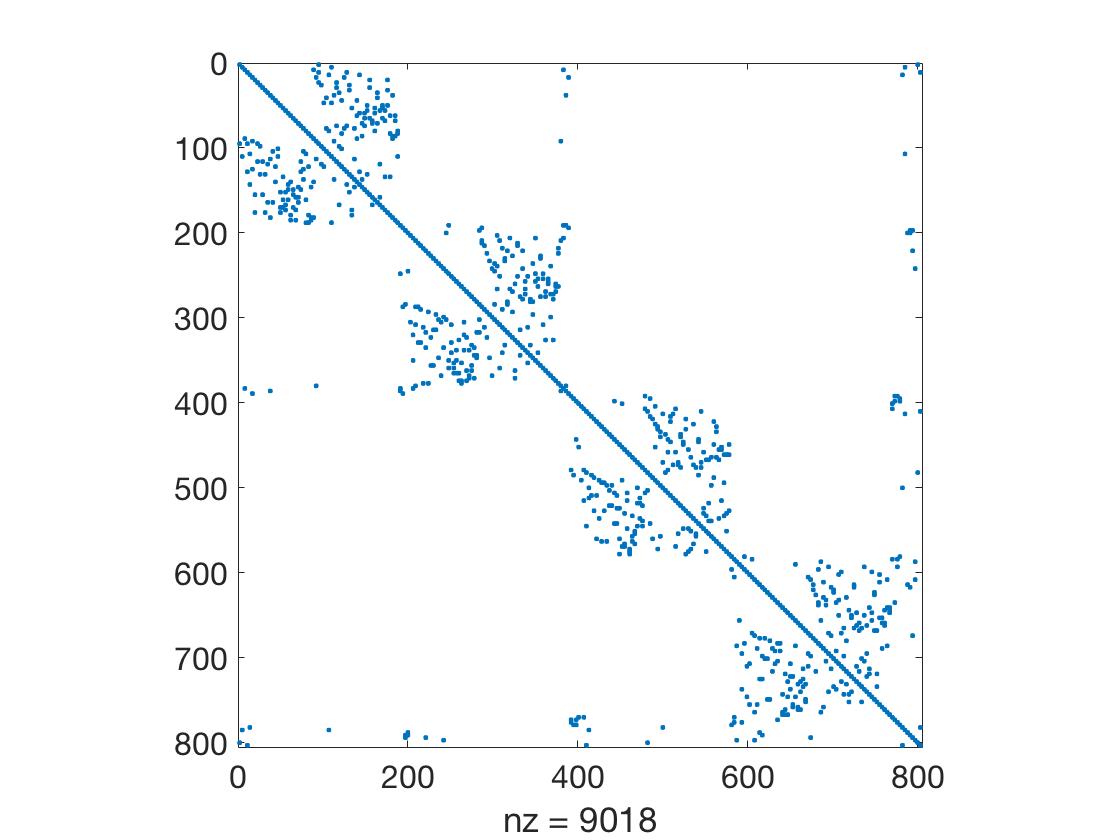}
\includegraphics[width = .5\textwidth]{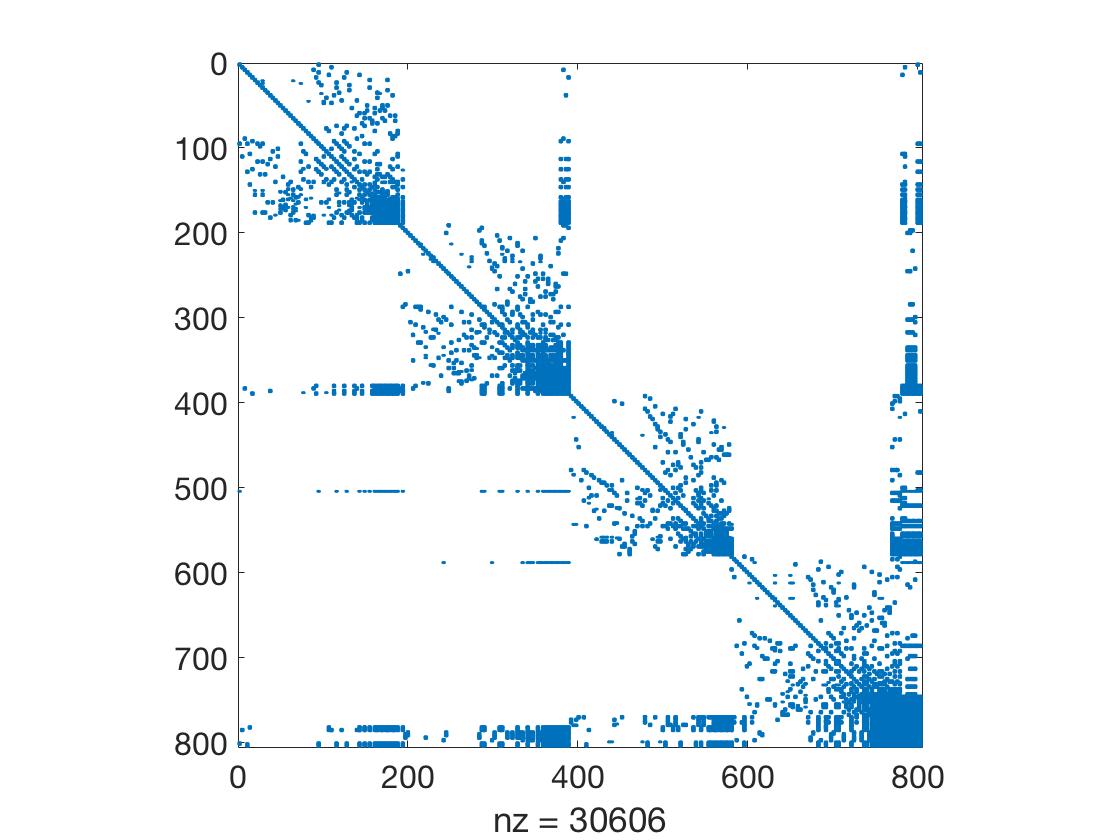}
\caption{\label{fig:Exp4ND5} (Left) sparsity structure of the global matrix 
after ND pivoting. (Right) sparsity structure of the combined LU decomposition 
the global matrix after ND pivoting.}
\end{figure}

\begin{figure}[!htb]
\includegraphics[width = .5\textwidth]{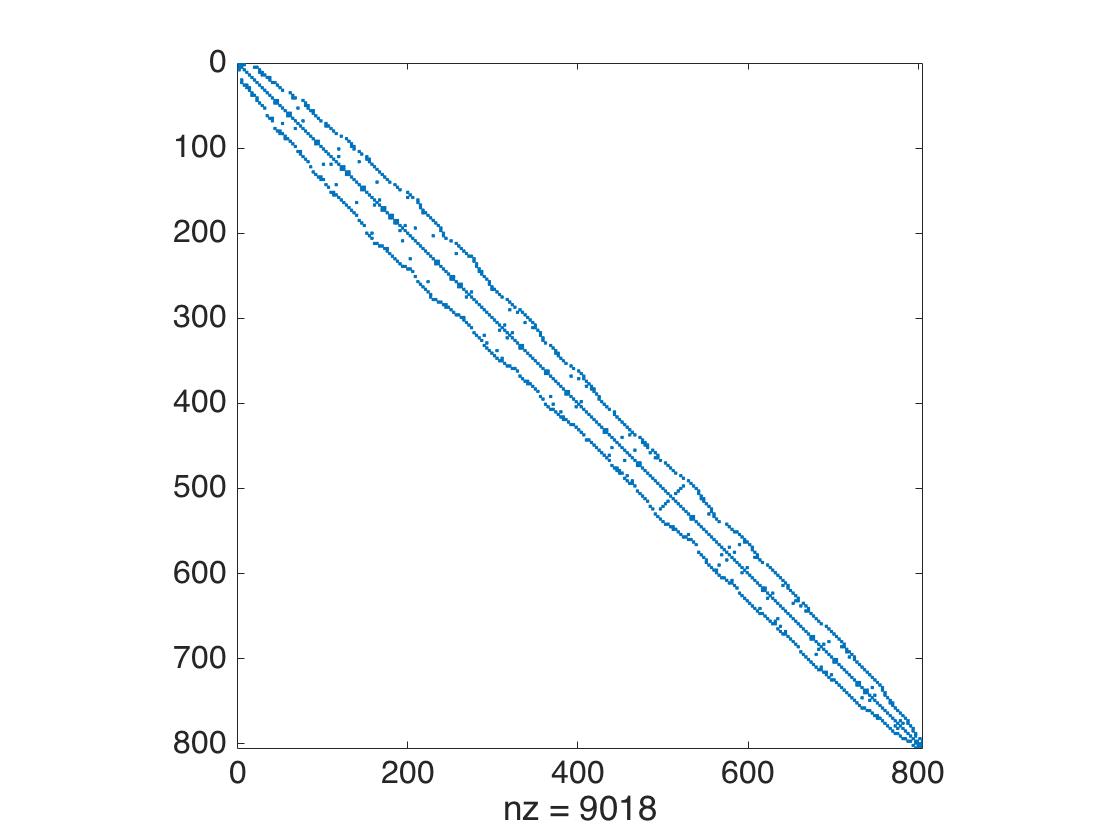}
\includegraphics[width = .5\textwidth]{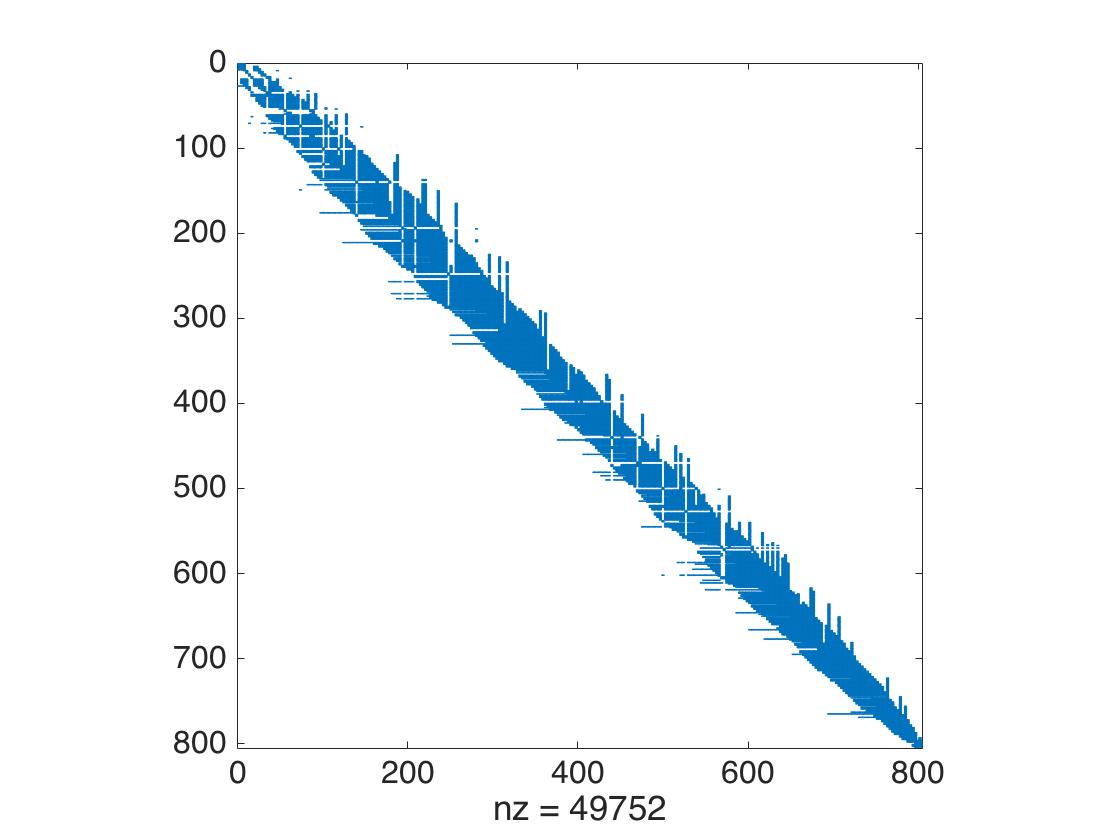}
\caption{\label{fig:Exp4RCM5} (Left) sparsity structure of the global matrix 
after RCM pivoting. (Right) sparsity structure of the combined LU 
decomposition the global matrix after RCM pivoting.}
\end{figure}

The sparsity structure of the matrix $A$ defined in \eqref{eq:LinSystem} along 
with its combined $LU$ factorization are given in Figure 
\ref{fig:Exp4NoPivot5}.  This is followed up with the sparsity structure of 
the pivoted matrix $PAP^T$ and the combined $LU$ factorization of this matrix 
for the AMD, ND, and RCM pivot strategies in Figure
\ref{fig:Exp4AMD5}--\ref{fig:Exp4RCM5}.  Tables \ref{tab:Exp4_a} and 
\ref{tab:Exp4_b} are used to summarize the data from the experiments in this 
section.  In particular, Table \ref{tab:Exp4_a} give the number of nonzero 
entries in the matrix $A$, its $LU$ factorization, and the $LU$ factorization 
of the matrix $PAP^T$ after AMD, ND, and RCM pivoting has been applied.  Table 
\ref{tab:Exp4_b} presents this nonzero data as a percent of the total number 
of entries in the matrix.

In contrast to experiments in Sections 
\ref{sec:NumericalExp1}--\ref{sec:NumericalExp3}, the RCM method does not 
perform well in the experiments in this section.  In fact, RCM increases the 
fill-in for the cases of $n = 5, 10, 15$.  RCM only reduces fill-in in the 
fine mesh case with $n = 20$, but even then the reduction is only by a factor 
of $1.13$. 

Similar to our previous experiments in Sections 
\ref{sec:NumericalExp1}--\ref{sec:NumericalExp3}, the AMD method reduces fill-
in at a greater rate than ND for a coarse mesh characterized by $n = 5$.  As 
the mesh is refined ND reduces fill-in at a rate greater than AMD.  For 
instance, when $n = 20$ ND reduces fill-in at a rate of $2.36$ while AMD 
reduces fill-in at a rate of $1.33$.


\begin{table}
\centering
\begin{tabular}{|c||c|c|c|}
\hline
n & Total Entries & Nonzero Entries & LU \\
\hline \hline
5 & $6.46 \times 10^5$ & 9018  & 49113 \\
\hline
10 & $8.50 \times 10^6$ & 33732 & $3.04 \times 10^5$ \\
\hline
15 & $4.45 \times 10^7$ & 78174 & $1.08 \times 10^6$ \\
\hline
20 & $1.44 \times 10^8$ & $1.42 \times 10^5$ & $3.69 \times 10^6$ \\
\hline
\end{tabular} \\[.5cm]
\begin{tabular}{|c||c|c|c|}
\hline
n & AMD & ND & RCM \\
\hline \hline
5 & 27702 & 30606 & 49752 \\
\hline
10 & $1.92 \times 10^5$ & $1.89 \times 10^5$ & $3.43 \times 10^5$ \\
\hline
15 & $8.73 \times 10^5$ & $5.90 \times 10^5$ & $1.29 \times 10^6$ \\
\hline
20 & $2.83 \times 10^6$ & $1.58 \times 10^6$ & $3.34 \times 10^6$ \\
\hline
\end{tabular}
\caption{(Above) Total entries, number of nonzero entries before finding the 
$LU$ decomposition, and the number of nonzero entries in the $LU$ 
decomposition.
(Below) Number of nonzero entries of the $LU$ decomposition after applying AMD, 
ND, and RCM pivoting.\label{tab:Exp4_a}}
\end{table}

\begin{table}
\centering
\begin{tabular}{|c||c|c|c|c|c|}
\hline
n & Nonzero Entries & LU & AMD & ND & RCM \\
\hline \hline
5 & 1.4 \% & 7.6 \% & 4.3 \% & 4.7 \% & 7.7 \% \\
\hline
10 & 0.4 \% & 3.6 \% & 2.3 \% & 2.2 \% & 4.03 \% \\
\hline
15 & 0.18 \% & 2.4 \% & 1.96 \% & 1.3 \% & 2.9 \% \\
\hline
20 & 0.098 \% & 2.6 \% & 1.96 \% & 1.1 \% & 2.3 \% \\
\hline
\end{tabular}
\caption{Percentage of nonzero entries before finding the $LU$ decomposition, 
in the $LU$ decomposition, and in the $LU$ decomposition after applying AMD, 
ND, and RCM pivoting.\label{tab:Exp4_b}}
\end{table}

\subsection{Numerical experiment 5} \label{sec:NumericalExp5}

In the final set of experiments of this paper, we use the same domain and mesh 
as Section \ref{sec:NumericalExp1}.  See Figure \ref{fig:Exp1Mesh} for 
examples of the computational mesh $\mathcal{T}_h$ used in this experiment.  
In contrast to the experiments in Section \ref{sec:NumericalExp1}, this 
section makes use of the finite dimensional function space $V^2_h$ in 
\eqref{eq:IPDG}, i.e. the IP-DG method used in this section involves piecewise 
quadratic solutions across the triangulation $\mathcal{T}_h$. This change in 
the approximation space changes the structure of the global matrix in 
\eqref{eq:LinSystem} and further leads to different results with regard to 
pivoting and fill-in.

\begin{figure}[!htb]
\includegraphics[width = .5\textwidth]{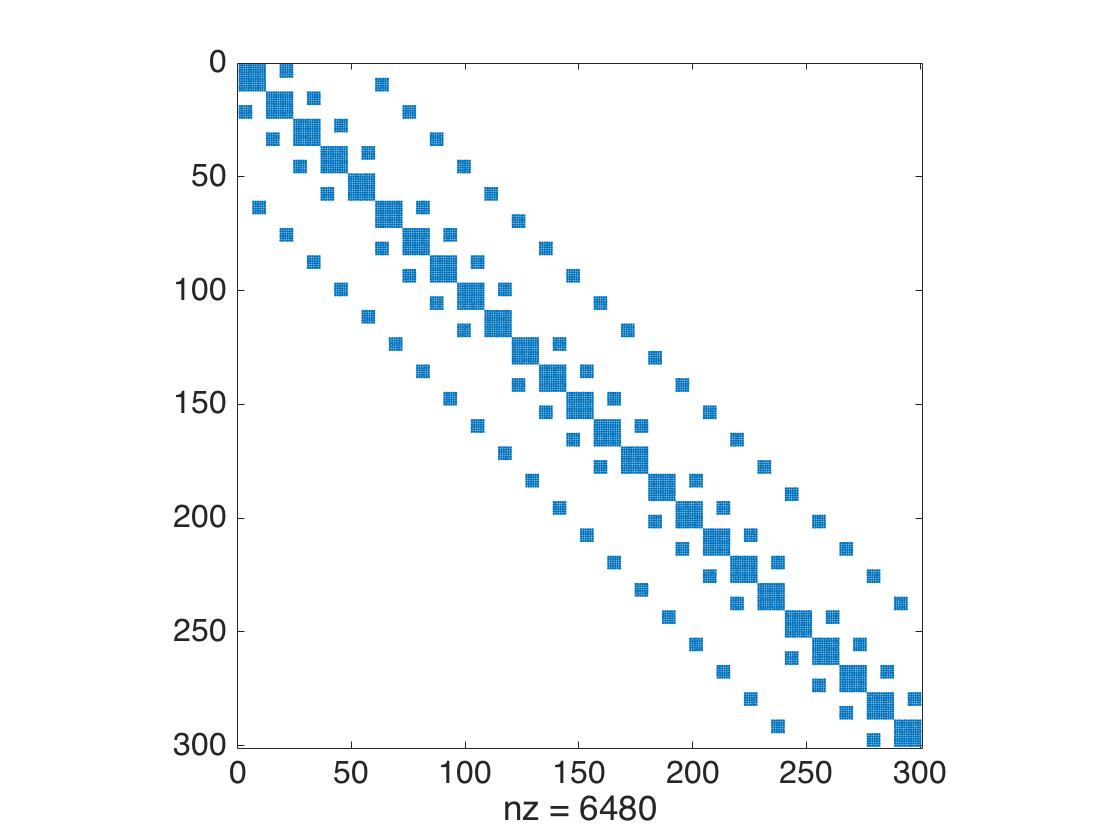}
\includegraphics[width = .5\textwidth]{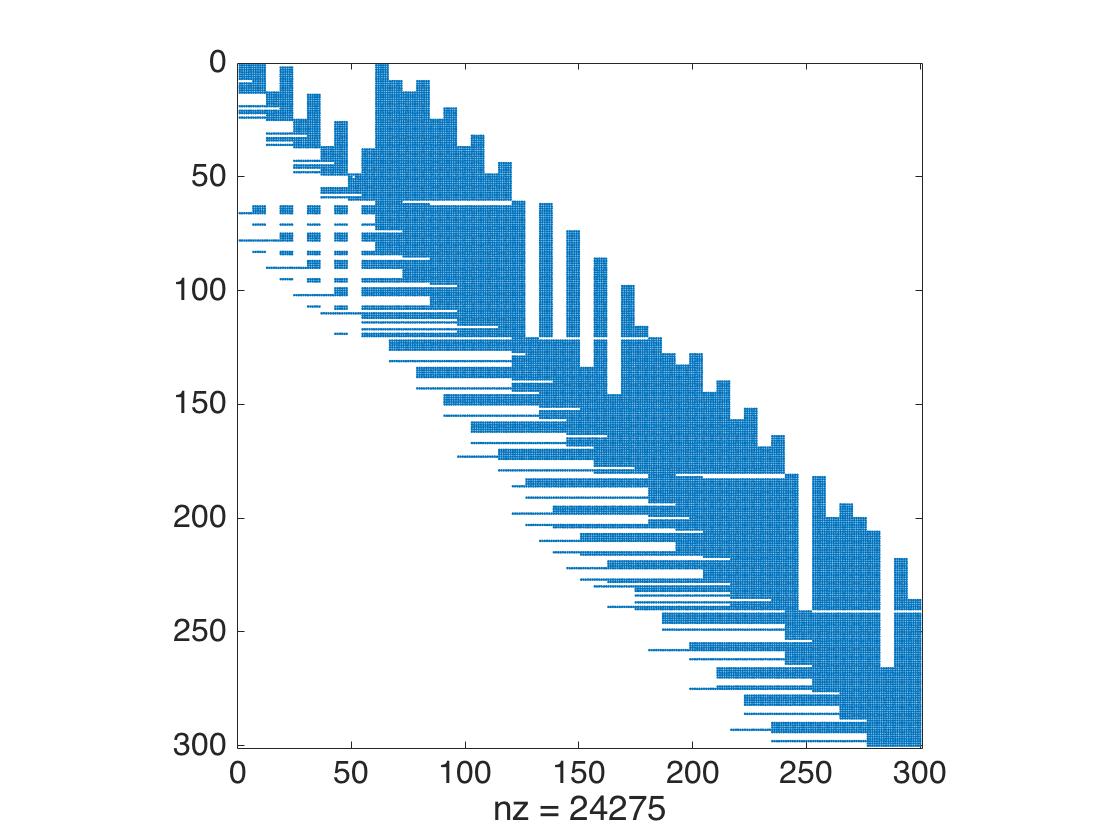}
\caption{\label{fig:Exp5NoPivot5} (Left) sparsity structure of the global 
matrix A produced by the IP-DG method with $n = 5$. (Right) sparsity structure 
of the combined LU decomposition of A.}
\end{figure}

\begin{figure}[!htb]
\includegraphics[width = .5\textwidth]{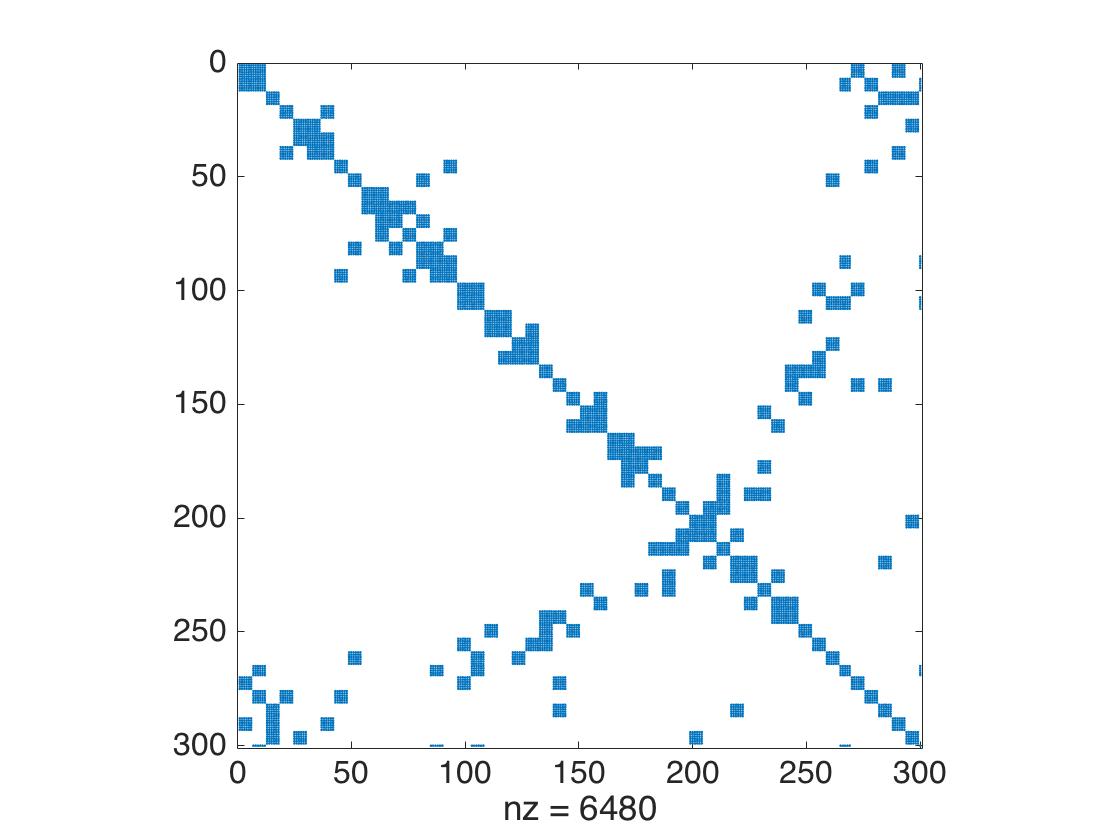}
\includegraphics[width = .5\textwidth]{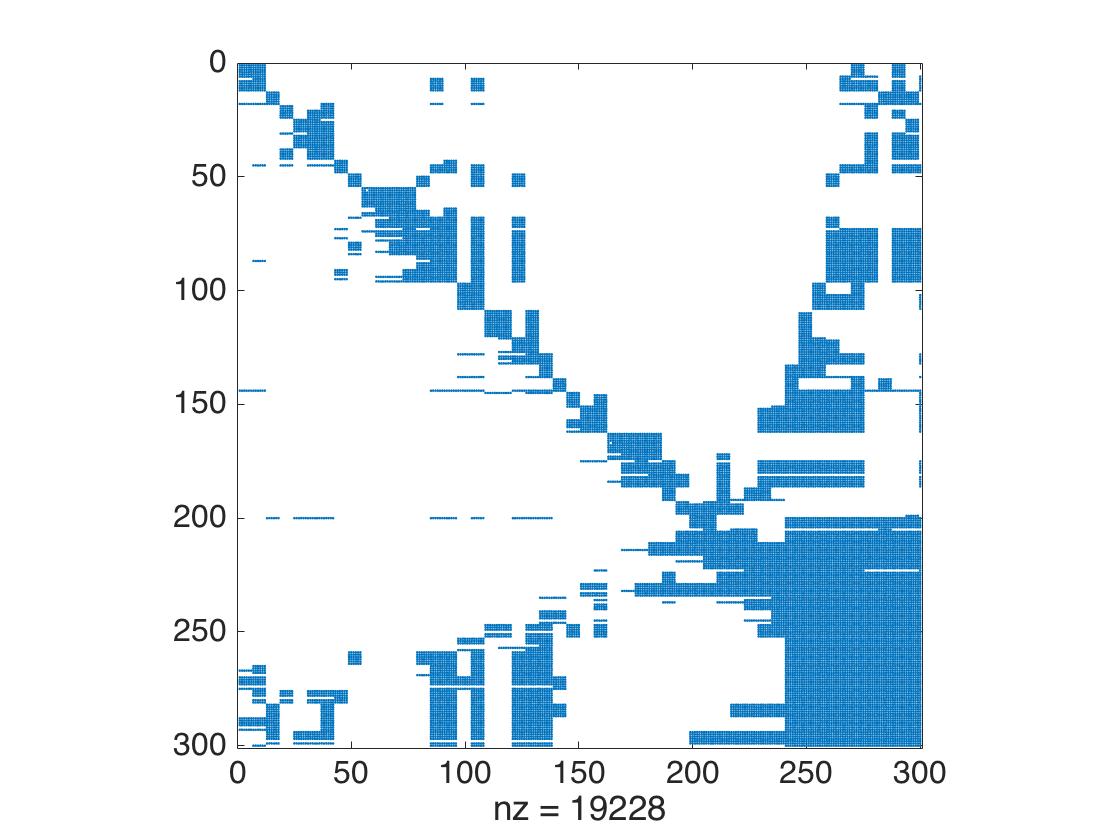}
\caption{\label{fig:Exp5AMD5} (Left) sparsity structure of the global matrix 
after AMD pivoting. (Right) sparsity structure of the combined LU 
decomposition the global matrix after AMD pivoting.}
\end{figure}

\begin{figure}[!htb]
\includegraphics[width = .5\textwidth]{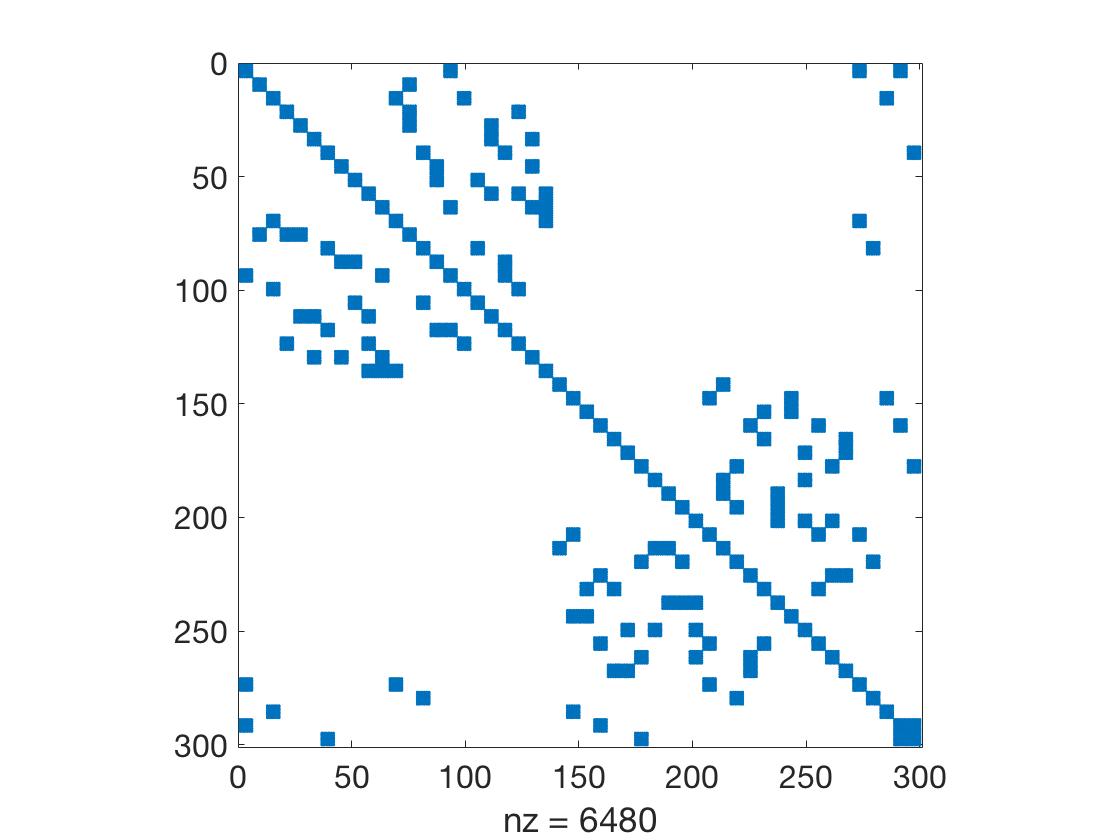}
\includegraphics[width = .5\textwidth]{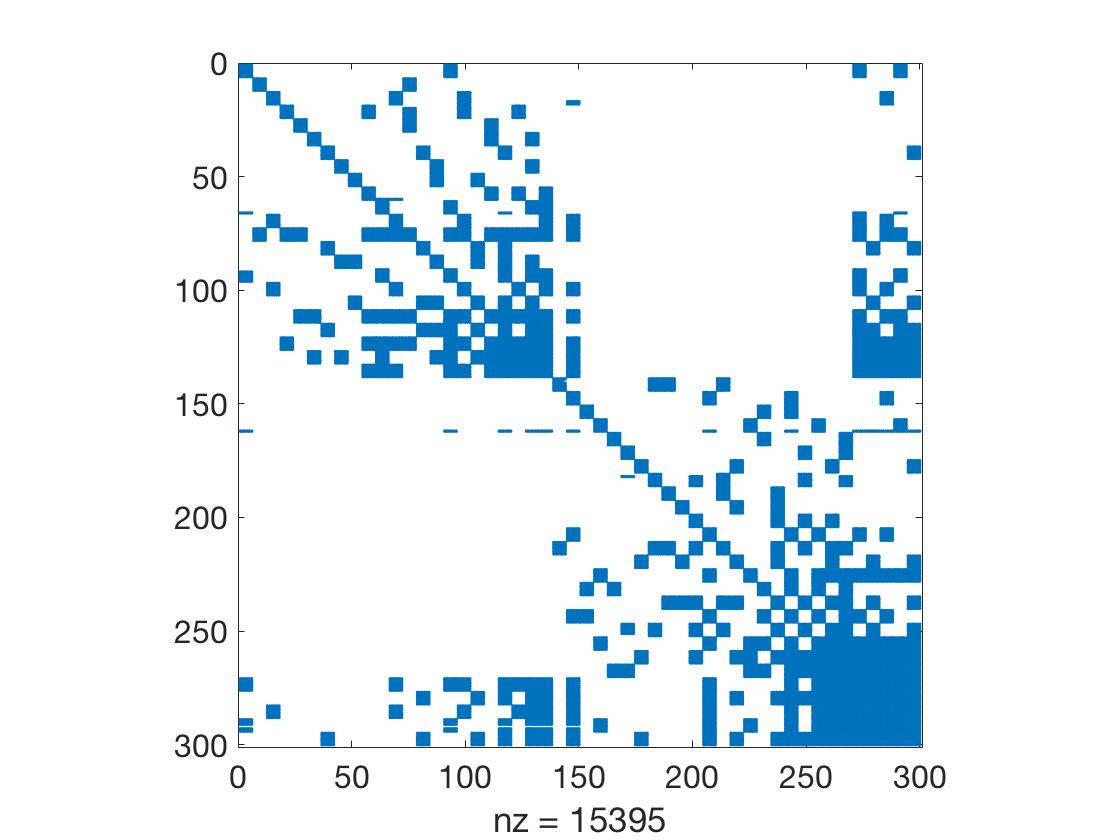}
\caption{\label{fig:Exp5ND5} (Left) sparsity structure of the global matrix 
after ND pivoting. (Right) sparsity structure of the combined LU decomposition 
the global matrix after ND pivoting.}
\end{figure}

\begin{figure}[!htb]
\includegraphics[width = .5\textwidth]{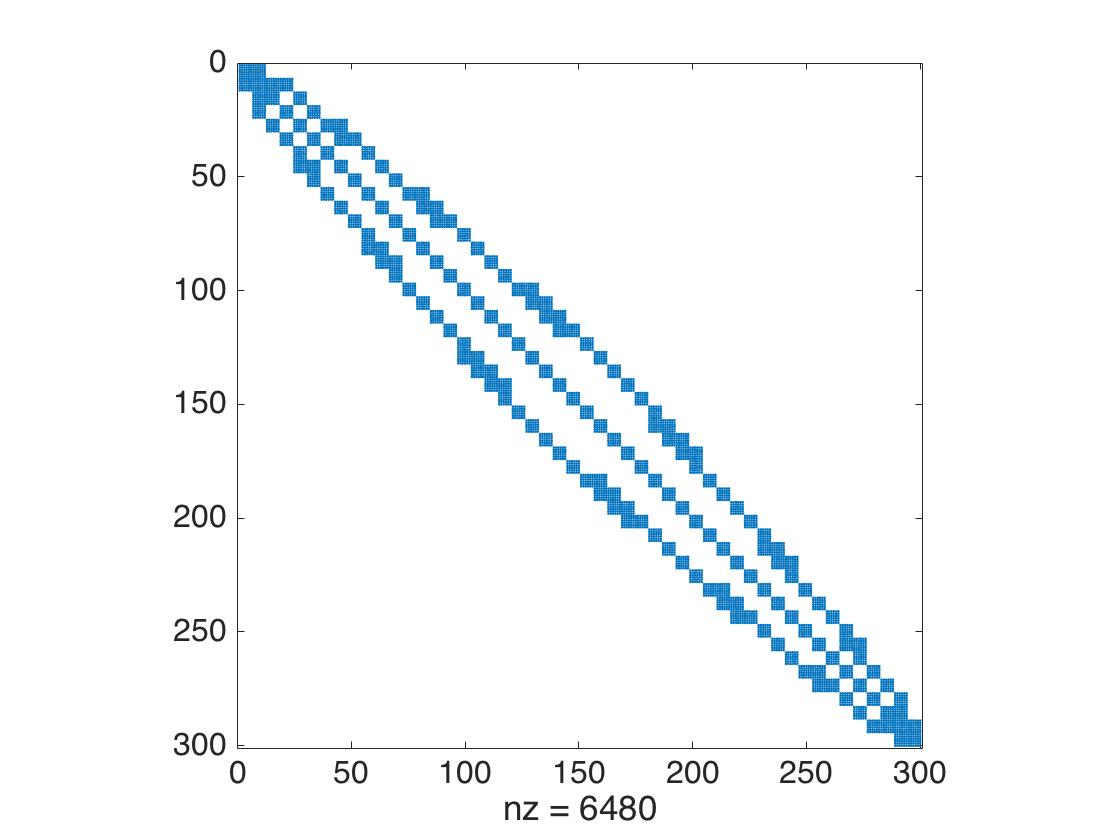}
\includegraphics[width = .5\textwidth]{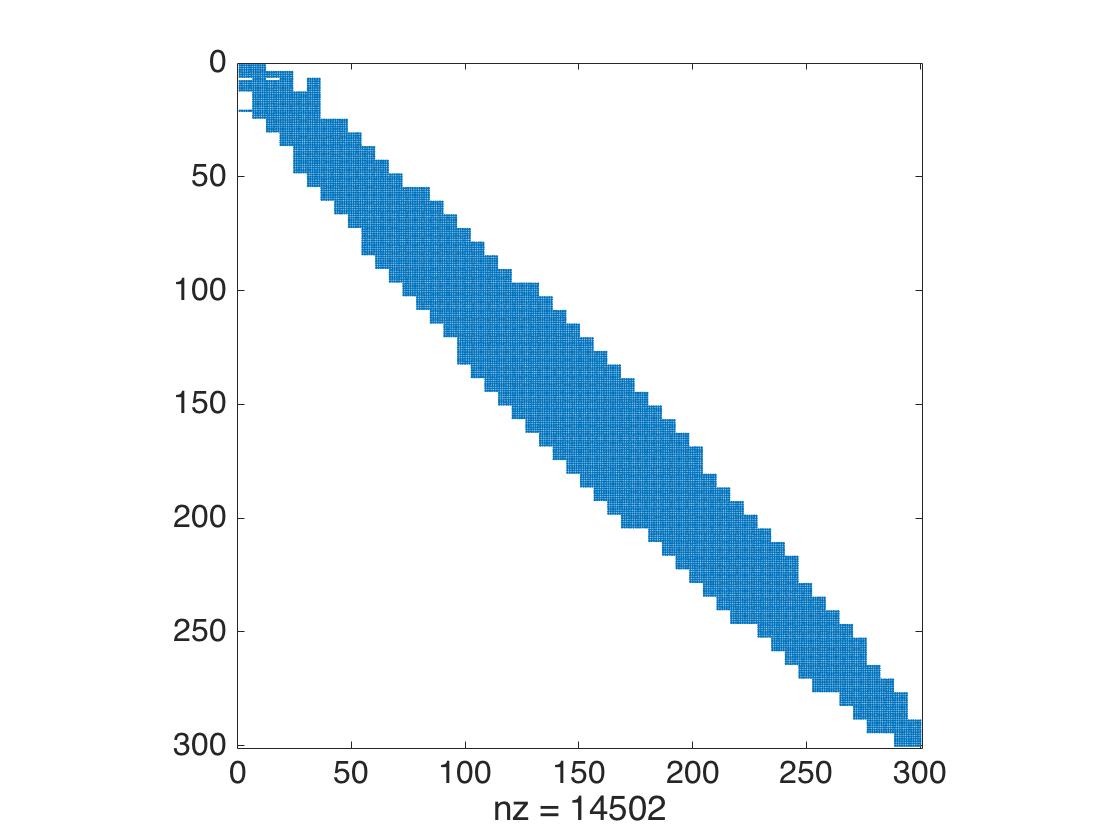}
\caption{\label{fig:Exp5RCM5} (Left) sparsity structure of the global matrix 
after RCM pivoting. (Right) sparsity structure of the combined LU 
decomposition the global matrix after RCM pivoting.}
\end{figure}

Figure \ref{fig:Exp5NoPivot5} shows the sparse structure of the matrix $A$ 
without pivoting and the sparse structure of the combined $LU$ factorization 
of $A$ when $n = 5$.  Figures \ref{fig:Exp5AMD5}--\ref{fig:Exp5RCM5} display 
the sparse structure of the pivoted matrix $PAP^T$ along with its combined 
$LU$ factorization for the AMD, ND, and RCM pivoting strategies when $n = 5$.  
From these figures we see differences in the fill-in produced in this section 
compared to those in Section \ref{sec:NumericalExp1}.  In particular, we see 
that the amount of fill-in from $LU$ factorization increased in this set of 
experiments.  Tables \ref{tab:Exp5_a} and \ref{tab:Exp5_b} summarize the 
results of the experiments using $n = 5, 10, 15, 20$ intervals along each side 
of the domain to generate the mesh $\mathcal{T}_h$.


\begin{table}[!htb]
\centering
\begin{tabular}{|c||c|c|c|}
\hline
n & Total Entries & Nonzero Entries & LU \\
\hline \hline
5 & 90000 & 6480 & 24275 \\
\hline
10 & $1.44 \times 10^6$ & 27360 & $2.00 \times 10^5$ \\
\hline
15 & $7.29 \times 10^6$ & 62640 & $7.09 \times 10^5$ \\
\hline
20 & $2.30 \times 10^7$ & $1.12 \times 10^5$ & $1.25 \times 10^6$ \\
\hline
\end{tabular} \\[.5cm]
\begin{tabular}{|c||c|c|c|}
\hline
n & AMD & ND & RCM \\
\hline \hline
5 & 19228 & 15395 & 14502 \\
\hline
10 & $3.26 \times 10^5$ & $2.18 \times 10^5$ & $1.16 \times 10^5$ \\
\hline
15 & $1.16 \times 10^6$ & $9.32 \times 10^5$ & $4.22 \times 10^5$ \\
\hline
20 & $3.02 \times 10^6$ & $1.74 \times 10^6$ & $9.91 \times 10^5$ \\
\hline
\end{tabular}
\caption{(Above) Total entries, number of nonzero entries before finding the 
$LU$ decomposition, and the number of nonzero entries in the $LU$ 
decomposition. (Below) Number of nonzero entries of the $LU$ decomposition 
after applying AMD, ND, and RCM pivoting.\label{tab:Exp5_a}}
\end{table}

\begin{table}[!htb]
\centering
\begin{tabular}{|c||c|c|c|c|c|}
\hline
n & Nonzero Entries & LU & AMD & ND & RCM \\
\hline \hline
5 & 7.2 \% & 26.97 \%  & 21.4 \%  & 17.1 \% & 16.1 \%  \\
\hline
10 & 1.9 \%  & 13.9 \%  & 22.7 \%  & 15.1 \% & 8.1 \%  \\
\hline
15 & 0.86 \%  & 9.7 \%  & 15.9 \%  & 12.8 \% & 5.8 \%  \\
\hline
20 & 0.49 \%  & 5.4 \%  & 13.1 \%  &  7.6 \% & 4.3 \%  \\
\hline
\end{tabular}
\caption{Percentage of nonzero entries before finding the $LU$ decomposition, 
in the $LU$ decomposition, and in the $LU$ decomposition after applying AMD, 
ND, and RCM pivoting.\label{tab:Exp5_b}}
\end{table}

Surprisingly, unlike previous experiments the RCM method out-performs both AMD 
and ND for all experiments in this section.  For a coarse mesh characterized 
by $n = 5$ all three methods reduce fill-in with a reduction rate of $1.26$ by 
AMD, $1.58$ by ND, and $1.674$ by RCM.  The performance of AMD and ND 
deteriorates as the mesh is refined in this set of experiments.  In fact, 
fill-in increases in the cases of $n = 10, 15, 20$ when using AMD and ND 
pivoting prior to the $LU$ factorization.  RCM continues to reduce fill-in 
when the mesh is refined, but for $n = 20$ the rate of reduction is decreased 
to $1.26$.

%
\section{Conclusions} \label{sec:Conclusions}
The results discussed in Section \ref{sec:NumericalExp} demonstrate that pivoting strategies are effective in reducing the fill-in due to LU factorization applied to the coefficient matrix generated by the IP-DG method applied to the Helmholtz problem \eqref{eq:PDE1}--\eqref{eq:PDE3}.  For experiments using a piecewise linear solution space $V^1_h$, AMD and ND reduced fill-in by similar amounts with AMD usually performing better than ND for coarse meshes and ND usually performing better than AMD for fine meshes (c.f. Sections \ref{sec:NumericalExp1}, \ref{sec:NumericalExp3}, and \ref{sec:NumericalExp4}).  In Section \ref{sec:NumericalExp2} ND performs better than AMD in both coarse and fine meshes.  For a piecewise linear solution space $V^1_h$ we observe that RCM does reduce fill-in, but not by as large of a factor as AMD and ND (c.f. Sections \ref{sec:NumericalExp1}--\ref{sec:NumericalExp4}).  Surprisingly, for a piecewise quadratic solution space $V^2_h$, RCM reduces fill-in by a larger factor than AMD and ND.  In fact, Section \ref{sec:NumericalExp5} demonstrates that AMD and ND increase fill-in for the solution space $V^2_h$ when using a fine mesh.

\bibliographystyle{plain}
\bibliography{bibliography}

\def\cprime{$'$}
\begin{thebibliography}{10}

\bibitem{Amestoy_Davis_Duff_96}
P.R. Amestoy, T.A. Davis, and I.S. Duff.
\newblock An approximate minimum degree ordering algorithm.
\newblock {\em SIAM J. Matrix Anal. Appl.}, 17:886 -- 905, 1996.

\bibitem{Amestoy_Davis_Duff_04}
P.R. Amestoy, T.A. Davis, and I.S. Duff.
\newblock Amd, an approximate minimum degree ordering algorithm.
\newblock {\em ACM Trans. Math. Software}, 30:381 -- 388, 2004.

\bibitem{Aziz_Kellogg_79}
A.K. Aziz and R.~B. Kellogg.
\newblock A scattering problem for the {Helmholtz} equation.
\newblock In {\em Advances in Computer Methods for Partial Differential
  Equations - III; Proceedings of the Third International Symposium, Bethlem,
  PA, June 20 - 22, 1979}, pages 93 -- 95, New Brunswick, NJ, 1979. IMACS.

\bibitem{Babuska_Sauter_97}
I.M. Babu\v{s}ka and S.A. Sauter.
\newblock Is the pollution effect of the {FEM} avoidable for the {Helmholtz}
  equation considering high wave numbers?
\newblock {\em SIAM J. Numer. Anal.}, 34(6):2392 -- 2423, 1997.

\bibitem{Buffa_Monk_08}
A.~Buffa and P.~Monk.
\newblock Error estimates for the ultra weak variational formulation of the
  {Helmholtz} equation.
\newblock {\em M2AN}, 42(6):925 -- 940, 2008.

\bibitem{Cessenat_Despres_98}
O.~Cessenat and B.~Despr\'{e}s.
\newblock Application of the ultra weak variational formulation to the two
  dimensional {Helmholtz} problem.
\newblock {\em SIAM J. Numer. Anal.}, 35:255 -- 299, 1998.

\bibitem{Cessenat_Despres_03}
O.~Cessenat and B.~Despr\'{e}s.
\newblock Using plane waves as base functions for solving time harmonic
  equations with the ultra weak variational formulation.
\newblock {\em J. Comput. Acoust.}, 11:227--238, 2003.

\bibitem{Chan_George_80}
W.M. Chan and A.~George.
\newblock A linear time implementation of the reverse cuthill-mckee algorithm.
\newblock {\em BIT}, 20:8 -- 14, 1980.

\bibitem{Cuthill_McKee_69}
E.~Cuthill and J.~McKee.
\newblock Reducing the bandwidth of sparse symmetric matrices.
\newblock In {\em Proceedings of the 24th Conference of the ACM}, pages 157 --
  172, Pinceton, NJ, 1969. Brandon Press.

\bibitem{Davis_2006}
T.~Davis.
\newblock {\em Direct Methods for Sparse Linear Systems}.
\newblock SIAM, Philadelphia, 2006.

\bibitem{Douglas_Santos_Sheen_Bennethum_93}
J.~Douglas~Jr., J.E. Santos, D.~Sheen, and L.S. Bennethum.
\newblock Frequency domain treatment of one-dimensional scalar waves.
\newblock {\em M$^3$AS}, 3(2):171 -- 194, 1993.

\bibitem{Douglas_Sheen_Santos_94}
J.~Douglas~Jr., D.~Sheen, and J.E. Santos.
\newblock Approximations of scalar waves in the space-frequency domain.
\newblock {\em M$^3$AS}, 4(4):509 -- 531, 1994.

\bibitem{Duff_Erisman_Reid_76}
I.S. Duff, A.M. Erisman, and J.K. Reid.
\newblock On {George's} nested dissection method.
\newblock {\em SIAM J. Numer. Anal.}, 13:686 -- 695, 1976.

\bibitem{Engquist_Majda_79}
B.~Engquist and A.~Majda.
\newblock Radiation bounday conditions for acoustic and elastics wave
  calculations.
\newblock {\em Comm. Pure Appl. Math.}, 32(3):314 -- 358, 1979.

\bibitem{Ernst_Gander_12}
O.~Ernst and M.~Gander.
\newblock Why it is difficult to solve {Helmholtz} problems with classical
  iterative methods?
\newblock In I.~Graham, T.~Hou, O.~Lakkis, and R.~Scheichl, editors, {\em
  Numerical Analysis of Multiscale Problems}, Lecture Notes in Computational
  Science and Engineering 83, pages 325 -- 363. Springer Verlag, 2012.

\bibitem{Feng_Lin_Lorton_15}
X.~Feng, J.~Lin, and C.~Lorton.
\newblock An efficient numerical method for acoustic wave scattering in random
  media.
\newblock {\em SIAM/ASA J. UQ}, 3:790 -- 822, 2015.

\bibitem{Feng_Wu_09}
X.~Feng and H.~Wu.
\newblock Discontinuous {Galerkin} methods for the {Helmholtz} equation with
  large wave numbers.
\newblock {\em SIAM J. Numer. Anal.}, 47:2872 -- 2896, 2009.

\bibitem{Feng_Wu_11}
X.~Feng and H.~Wu.
\newblock hp-discontinuous {Galerkin} methods for the {Helmholtz} equation with
  large wave numbers.
\newblock {\em Math. Comp.}, 80:1997 -- 2024, 2011.

\bibitem{George_73}
A.~George.
\newblock Nested dissection of a regular finite element mesh.
\newblock {\em SIAM. J. Numer. Anal.}, 10:345 -- 363, 1973.

\bibitem{George_Ng_85}
A.~George and E.~Ng.
\newblock An implementation of gaussian elimination with partial pivoting for
  sparse systems.
\newblock {\em SIAM. J. Sci. Statist. Comput.}, 6:390 -- 409, 1985.

\bibitem{Gilbert_94}
J.R. Gilbert.
\newblock Predicting structure in sparse matrix computations.
\newblock {\em SIAM J. Matrix Anal. Appl.}, 15:62 -- 79, 15.

\bibitem{Gittelson_Hiptmair_Perugia_09}
C.~Gittelson, R.~Hiptmair, and I.~Perugia.
\newblock Plane wave discontinuous {Galerkin} methods: Analysis of the
  h-version.
\newblock {\em M2AN}, 43:297 -- 331, 2009.

\bibitem{Hecht_12}
F.~Hecht.
\newblock New development in freefem++.
\newblock {\em J. Numer. Math.}, 20(3-4):251--265, 2012.

\bibitem{Hiptmair_Moiola_Perugia_11_2}
R.~Hiptmair, A.~Moiola, and I.~Perugia.
\newblock Plane wave discontinuous {Galerkin} methods for the 2d {Helmholtz}
  equation: Analysis of the p-version.
\newblock {\em SIAM J. Numer. Anal.}, 49:264 -- 284, 2011.

\bibitem{Huttunen_Monk_07}
T.~Huttunen and P.~Monk.
\newblock The use of plane waves to approximate wave propagation in anisotropic
  media.
\newblock {\em J. Comput. Math.}, 25:350 -- 367, 2007.

\bibitem{Huttunen_Monk_Collino_Kaipio_04}
T.~Huttunen, P.~Monk, F.~Collino, and J.P. Kaipio.
\newblock The ultra-weak variational formulation for elastic wave problems.
\newblock {\em SIAM J. Sci. Comput.}, 25(5):1717 -- 1742, 2004.

\bibitem{Ihlenburg_Babuska_95}
F.~Ihlenburg and I.~Bab\v{u}ska.
\newblock Finite element solution of the {Helmholtz} equation with high wave
  number. the h-version of {FEM}.
\newblock {\em Comp. Math. Appl.}, 30(9):9 -- 37, 1995.

\bibitem{Liu_Sherman_76}
W.-H Liu and A.H. Sherman.
\newblock Comparitive analysis of the {Cuthill-McKee} and the reverse {Cuthill
  McKee} ordering algorithms for sparse matrices.
\newblock {\em SIAM J. Numer. Anal.}, 13:198 -- 213, 1976.

\bibitem{Luostari_Huttunen_Monk_09}
T.~Luostari, T~Huttunen, and P.~Monk.
\newblock Plane wave methods for approximating the time harmonic wave equation.
\newblock In B.~Engquist, A.~Fokas, E.~Hairer, and A.~Iserles, editors, {\em
  Highly oscillatory problems}, volume 366 of {\em London Mathematcal Society
  Lecture Notes}, pages 127 -- 153. Cambridge University Press, Cambridge,
  2009.

\bibitem{Markowitz_57}
H.M. Markowitz.
\newblock The elimination form of the inverse and its application to linear
  programming.
\newblock {\em Management Sci.}, 3:255 -- 269, 1957.

\bibitem{Matlab_AMD}
MathWorks.
\newblock {Matlab Documentation}: Approximate minimum degree permutation.
\newblock \url{https://www.mathworks.com/help/matlab/ref/amd.html}, 1994 -
  2019.
\newblock Accessed: 2019-01-31.

\bibitem{Matlab_ND}
MathWorks.
\newblock {Matlab Documentation}: Nested dissection permutation.
\newblock \url{https://www.mathworks.com/help/matlab/ref/dissect.html}, 1994 -
  2019.
\newblock Accessed: 2019-01-31.

\bibitem{Matlab_RCM}
MathWorks.
\newblock {Matlab Documentation}: Sparse reverse cuthill-mckee ordering.
\newblock \url{https://www.mathworks.com/help/matlab/ref/symrcm.html}, 1994 -
  2019.
\newblock Accessed: 2019-01-31.

\bibitem{Melenk_Sauter_11}
J.M. Melenk and S.~Sauter.
\newblock Wavenumber explicit convergence analysis for {Galerkin}
  discretization of the {Helmholtz} equation.
\newblock {\em SIAM J. Numer. Anal.}, 49(3):1210 -- 1243, 2011.

\bibitem{Moiola_09}
A.~Moiola.
\newblock Approximation properties of plan wave spaces and application to the
  analysis of the plane wave discontinuous {Galerkin} method.
\newblock In {\em Technical report 2009-06, Seminar f\"{u}r Angewandte
  Mathematik}. ETH Z\"{u}rich, 2009.

\bibitem{Riviere_08}
B.~Rivi\`{e}re.
\newblock {\em Discontinuous {Galerkin} methods for solving elliptic and
  parabolic equations: theory and implementation}.
\newblock SIAM, Philadelphia, PA, 2008.

\bibitem{Wu_13}
H.~Wu.
\newblock Pre-asymptotic error analysis of {CIP}-{FEM} and {FEM} for the
  {Helmholtz} equation with high wave number. part i: linear version.
\newblock {\em IMA J. Numer. Anal.}, 34(3):1266 -- 1288, 2013.

\end{thebibliography}

\end{document}